\documentclass[paper,10pt]{IEEEtran}
\usepackage[T1]{fontenc}
\usepackage{amsmath,amssymb,amsfonts}
\usepackage{graphicx}
\usepackage{xcolor}
\usepackage{textcomp}
\usepackage{algorithm}
\usepackage{algpseudocode}
\usepackage{array}
\usepackage{tabularx}
\usepackage[acronym]{glossaries}
\newacronym{3d}{3-D}{three-dimensional}
\newacronym{3gpp}{3GPP}{3rd Generation Partnership Project}
\newacronym{5g}{5G}{fifth generation}
\newacronym{6g}{6G}{sixth generation}
\newacronym{ao}{AO}{alternating optimization}
\newacronym{aoa}{AoA}{angle of arrival}
\newacronym{aod}{AoD}{angle of departure}
\newacronym{as}{AS}{antenna selection}
\newacronym{awgn}{AWGN}{additive white gaussian noise}
\newacronym{b5g}{B5G}{beyond fifth generation}
\newacronym[plural = BSs, firstplural = base statios (BSs)]{bs}{BS}{base station}
\newacronym{ber}{BER}{bit-error-rate}
\newacronym{bfs}{BFS}{backward-forward sweep}
\newacronym{csi}{CSI}{channel state information}
\newacronym{chest}{CHEST}{channel estimation}
\newacronym{ccdf}{CCDF}{complementary cumulative distribution function}
\newacronym{ccm}{CCM}{complex circle manifold}
\newacronym{cs}{CS}{channel sounding}
\newacronym{ccp}{CCP}{convex-concave procedure}
\newacronym{dl}{DL}{downlink}
\newacronym{dcp}{DCP}{disciplined convex programming}
\newacronym{ee}{EE}{energy efficiency}
\newacronym{em}{EM}{electromagnetic}
\newacronym{er}{ER}{ergodic rate}
\newacronym{epds}{EPDS}{electrical power distribution systems}
\newacronym{epts}{EPTS}{electrical power transmission systems}
\newacronym{fbs}{FBS}{forward-backward sweep}
\newacronym{fp}{FP}{fractional programming}
\newacronym{fdm}{FDM}{finite difference method}
\newacronym{ga}{GA}{genetic algorithm}
\newacronym{gf}{GF}{grant-free}
\newacronym{gwo}{GWO}{grey wolf optimization}
\newacronym{hbf}{HBF}{hybrid beamforming}
\newacronym[plural = HEMs]{hem}{HEM}{{\it heuristic evolutionary}}
\newacronym{iao}{i-AO}{iterative-alternating optimization}
\newacronym{iid}{i.i.d.}{independent and identically distributed}
\newacronym{iot}{IoT}{internet of things}
\newacronym{isomap}{Isomap}{Isometric Mapping}
\newacronym{kkt}{KKT}{Karush–Kuhn–Tucker}
\newacronym[plural = KPIs, firstplural = key performance indicators (KPIs)]{kpi}{KPI}{key performance indicator}
\newacronym{ldt}{LDT}{Lagrangian Dual Transform}
\newacronym{los}{LoS}{line-of-sight}
\newacronym{lle}{LLE}{Locally Linear Embedding}

\newacronym{mimo}{MIMO}{multiple-input multiple-output}
\newacronym{m-mimo}{mMIMO}{massive MIMO}
\newacronym{mmse}{MMSE}{minimum mean squared error}
\newacronym{mrc}{MRC}{maximum ratio combiner}
\newacronym{mcs}{MCs}{Monte-Carlo Simulation}
\newacronym{mrt}{MRT}{maximum ratio transmission}
\newacronym{mu}{MU}{multi-user}
\newacronym{mo}{MO}{Manifold optimization}
\newacronym{mds}{MDS}{Multidimensional Scaling }
\newacronym{mmtc}{mMTC}{massive machine type communications}
\newacronym{mm}{MM}{majorization minimization}
\newacronym{mtc}{MTC}{machine type communications}
\newacronym{ml}{ML}{machine-learning}
\newacronym{mip}{MIP}{mixed integer programming}
\newacronym{milp}{MILP}{mixed integer linear programming}
\newacronym{minlp}{MINLP}{mixed integer nonlinear programming}
\newacronym{np-complete}{NP-complete}{nondeterministic polynomial-time complete}
\newacronym{noma}{NOMA}{non-orthogonal multiple access}
\newacronym{nlos}{NLoS}{non line-of-sight}
\newacronym[plural=NNs]{nn}{NN}{neural network}
\newacronym{of}{OF}{objective function}
\newacronym{olp}{OLP}{optimal linear precoder}
\newacronym{pa}{PA}{power amplifier}
\newacronym{pc}{PC}{pilot contamination}
\newacronym{pbo}{PBO}{passive beamforming optimization}
\newacronym{pbf}{PBF}{passive beamforming}
\newacronym{pca}{PCA}{principal component analysis}
\newacronym{prf}{PRF}{pilot-reuse factor}
\newacronym{pso}{PSO}{particle swarm optimization}
\newacronym{qbbfs}{QBBFS}{quadratic-based backward-forward sweep}
\newacronym{qos}{QoS}{Quality of Service}
\newacronym{rf}{RF}{radio frequency}
\newacronym{re}{RE}{resource efficiency}
\newacronym{rm}{RM}{Riemannian manifold}
\newacronym[plural=RISs]{ris}{RIS}{reconfigurable intelligent surface}
\newacronym{rcg}{RCG}{Riemann conjugate gradient}
\newacronym{ra}{RA}{random access}
\newacronym{rn}{RN}{Riemannian Newton}
\newacronym{rtr}{RTR}{Riemannian Trust-Region}
\newacronym{rgd}{RGD}{Riemannian Gradient Descent}
\newacronym{se}{SE}{spectral efficiency}
\newacronym{sfp}{SFP}{sequential fractional programming}
\newacronym{sinr}{SINR}{signal-to-interference-plus-noise ratio}
\newacronym{snr}{SNR}{signal-to-noise ratio}
\newacronym{sr}{SR}{sum rate}
\newacronym{sdr}{SDR}{semidefinite relaxation}
\newacronym{sca}{SCA}{sucessive convex approximation}
\newacronym{socp}{SOCP}{second-order cone programming}
\newacronym{sll}{SLL}{side lobe level}
\newacronym{sdp}{SDP}{semidefinite programming}
\newacronym{sg}{SG}{sub-gradient}
\newacronym{sa}{SA}{simulated annealing }
\newacronym{tdd}{TDD}{time-division duplex}
\newacronym{t-sne}{t-SNE}{t-Distributed Stochastic Neighbor Embedding}
\newacronym[plural=UEs, firstplural=users' equipment (UEs)]{ue}{UE}{user's equipment}
\newacronym{ul}{UL}{uplink}
\newacronym{ula}{ULA}{uniform linear array}
\newacronym{upa}{UPA}{uniform planar array}
\newacronym{uspa}{USPA}{uniform squared planar array}
\newacronym{ut}{UT}{unit terminal}
\newacronym{umap}{UMAP}{Uniform Manifold Approximation and Projection}
\newacronym{vr}{VR}{visibility region}
\newacronym{wbd}{WBD}{Wide-beam design}
\newacronym{xl-mimo}{XL-MIMO}{extra-large scale massive MIMO}
\newacronym{zf}{ZF}{Zero-Forcing}
\newacronym{zft}{ZFT}{zero-forcing transmission}

\usepackage[colorlinks=true, linkcolor=blue, citecolor=blue, urlcolor=blue]{hyperref}
\usepackage{subcaption}
\usepackage{siunitx}
\usepackage{soul}
\usepackage{comment}

\usepackage{ragged2e} 

\definecolor{gray}{rgb}{0.5,0.5,0.5}

\definecolor{teal}{RGB}{0, 128, 128}  
\newcommand{\colo}{\textcolor{teal}}

\newcommand{\dm} {\, \colo{$\diamond$}\, }
\renewcommand{\algorithmiccomment}[1]{{\color{gray} // #1}}

\newcolumntype{P}[1]{>{\centering\arraybackslash}p{#1}}
\newcolumntype{M}[1]{>{\centering\arraybackslash}m{#1}}

\definecolor{teal}{RGB}{0, 128, 128}
\definecolor{gray}{rgb}{0.5,0.5,0.5}

\renewcommand{\algorithmiccomment}[1]{{\color{gray} // #1}}

\def\BibTeX{{\rm B\kern-.05em{\sc i\kern-.025em b}\kern-.08em
T\kern-.1667em\lower.7ex\hbox{E}\kern-.125emX}}

\begin{document}

\title{Manifolds in Power Systems Optimization}

\author{
\IEEEauthorblockN{Lucca Rodrigues Pinto\IEEEauthorrefmark{1},
Wilson de Souza Junior\IEEEauthorrefmark{1}, 
Jaime Laelson Jacob\IEEEauthorrefmark{1}, 
Luis Alfonso Gallego Pareja\IEEEauthorrefmark{1}, 
Taufik Abrão\IEEEauthorrefmark{1}}
\IEEEauthorblockA{\IEEEauthorrefmark{1}State University of Londrina (UEL), CEP: 86057-970, Londrina, PR, Brazil \\
Email: taufik@uel.br}
}

\maketitle

\begin{abstract}
\gls{mo} is a powerful mathematical framework that can be applied to solving complex optimization problems with \glspl{of} and constraints on complex geometric structures, which is particularly useful in advanced power systems. We explore the application of \gls{mo} techniques, which offer a robust framework for solving complex, non-convex optimization problems in \gls{epds} and \gls{epts}, particularly for power flow analysis. This paper introduces the principles of \gls{mo} and demonstrates its advantages over conventional methods by applying it to power flow optimization. For \gls{epds}, a cost function derived from a \gls{bfs} algorithm is optimized using the Manopt toolbox, yielding high accuracy and competitive computational times on 14-bus, 33-bus, and 69-bus systems when compared to established solvers. Similarly, for \gls{epts}, \gls{mo} applied via Manopt to 3-bus and 4-bus systems effectively solves power flow equations, matching traditional methods such as Newton-Raphson in performance. The study highlights that tools such as Manopt can mitigate implementation complexities, positioning MO as an efficient and accessible tool for power system analysis and potentially broader planning applications. The paper provides a comprehensive tutorial on \gls{mo}, detailing its theoretical foundations, practical methodologies, and specific applications in power systems, particularly in power flow optimization.
\end{abstract}


\maketitle


\section{INTRODUCTION}

\IEEEPARstart{T}{he} core strength of \gls{mo} methods lies in their superior ability to tackle complex optimization problems, outperforming conventional methods in efficiency and effectiveness, particularly for the intricate non-convex optimization challenges inherent in \gls{epds} and \gls{epts} with their specific power flow considerations.




Many optimization problems in power systems are characterized by non-convex constraints, which frequently arise from the physical laws governing power flow and the operational limits of the network. Traditional optimization methods can find these constraints challenging to navigate efficiently. \gls{mo} offers a more inherently suitable approach by treating the set of feasible solutions—defined by these very constraints—as a smooth manifold. This paradigm allows \gls{mo} to directly leverage the problem's intrinsic geometric properties, leading to more efficient and often more robust solution methods.

A key technique enabling this is local linearization. While a manifold can be complex globally, it is locally Euclidean. This means that in the immediate vicinity of any point on the manifold, the space behaves much like standard Euclidean space. This property is crucial because it allows for the adaptation and generalization of powerful iterative optimization algorithms, such as gradient descent and Newton's method, to these otherwise curved and constrained spaces. The methodologies applied in this paper for power flow optimization indeed capitalize on such generalized gradient-based approaches.

Furthermore, power system models often involve a large number of variables and interdependent equations, creating a complex, high-dimensional landscape for optimization. \gls{mo} provides tools and perspectives to manage this complexity. Finally, and central to its application in this paper, \gls{mo} methods offer profound versatility in integrating the specific, often intricate, constraints of a given problem directly into the optimization process. Instead of treating power flow equations as separate, challenging boundary conditions, \gls{mo} allows these to define the very fabric of the manifold on which the optimization occurs. This adaptability in conforming to the unique structure of power system problems is a significant advantage, enabling the tailored and effective optimization strategies demonstrated herein.

Common approaches in the literature to tackle the complexities of the power flow model often involve linearization, which includes discretizing nonlinear terms \cite{FP:ref6}. However, such linearization introduces a dependency on parameters that must be tuned using the original nonlinear, non-convex model, which can compromise solution accuracy and introduce uncertainties. Another strategy, convexification of the power flow problem into models like \gls{socp}, offers benefits for planning and expansion, as combinatorial strategies like network reconfiguration can inherit this convexity \cite{convex1}, \cite{convex2}. Despite these advancements, both linearized and convex models encounter significant hurdles regarding solver accessibility. Solvers that support the integer variables standard in these strategies are often commercial. This reliance on often expensive commercial software for robust solutions presents a notable challenge, as freely available global nonlinear solvers for integer variables are scarce, and local solvers typically only ensure first-order KKT conditions, which may not be sufficient.

Given these challenges with traditional modeling and solver availability, alternative approaches are warranted. Direct application of nonlinear models, especially in the context of \gls{mo}, can face difficulties with initial point selection and computational demands, particularly as power flow calculations involve interdependent variables across shared buses, necessitating iterative computations akin to a \gls{bfs} approach \cite{Zimmerman1996}. The \gls{bfs} algorithm itself is a well-established and computationally efficient method for power flow analysis in \gls{epds}, especially for radial or weakly meshed systems, often outperforming traditional Newton-Raphson or Zbus methods in speed and robustness. For instance, studies have shown BFS, particularly with enhancements, can be significantly faster, and its variants like \gls{fbs} and \gls{qbbfs} have demonstrated superior performance in scenarios with high resistance-to-reactance ratios or multi-phase systems \cite{Muruganantham2016}, \cite{Mahmoud2016}. While a \gls{bfs} algorithm offers a robust alternative to direct modeling for power flow calculations, integrating its procedural nature into an \gls{mo} framework requires the thoughtful derivation of an appropriate cost function, a critical step explored in this paper.

\subsection{MOTIVATION AND CONTRIBUTION}

The optimization of modern \gls{epds} and \gls{epts} frequently involves solving complex, non-convex problems, particularly in power flow analysis. Traditional optimization methods often struggle with the inherent non-linearities and constraints of these systems, leading to computationally intensive processes, reliance on linearization or convex relaxation techniques that may sacrifice accuracy, or the use of heuristic methods that lack convergence guarantees. While \gls{mo} offers a theoretically powerful alternative by recasting constrained problems as unconstrained ones on smooth manifolds, its practical implementation can appear daunting due to the specialized knowledge of differential geometry typically required.

This is where toolboxes like Manopt become pivotal. A primary motivation for this work stems from the accessibility and flexibility Manopt brings to \gls{mo}. Manopt significantly lowers the barrier to applying sophisticated \gls{mo} techniques by requiring the user to primarily focus on defining a well-structured problem by stating a cost function and its well-defined gradient. Furthermore, Manopt operates within familiar environments like MATLAB or Python, meaning users are not confined by a restrictive, specialized optimization syntax. Instead, they can leverage the full expressive power and extensive libraries of these environments to formulate their cost functions, gradients, and any auxiliary computations. This flexibility motivates the exploration and practical application of \gls{mo} to complex power system problems that might otherwise be considered too intricate to implement from first principles. Building on this motivation, the contributions of this paper are multifaceted:
\begin{itemize}
\item[\dm] Practical Application of \gls{mo} via Manopt: We demonstrate the successful and practical application of \gls{mo}, facilitated by the Manopt toolbox, to solve non-convex power flow optimization problems in both \gls{epds} and \gls{epts}. This includes novel problem formulations where:
\begin{itemize}
\item[\dm] For \gls{epds}, a cost function is derived from the \gls{bfs} algorithm, with constraints handled via penalty methods, showcasing a flexible integration of procedural algorithms within the \gls{mo} framework.

\item[\dm] For \gls{epts}, standard power flow equations are directly used to define equality constraints, forming a cost function suitable for \gls{mo}.
\end{itemize}
\item[\dm] Performance Validation and Benchmarking: Through detailed case studies on standard test systems (14-bus, 33-bus, and 69-bus for \gls{epds}; and 3-bus and 4-bus systems for \gls{epts}), we rigorously validate the proposed \gls{mo} approaches. We show that \gls{mo}, implemented in Manopt, achieves high accuracy and exhibits computational performance that is highly competitive with, and in several instances comparable or superior to, well-established commercial solvers (like Knitro, CPLEX, Gurobi) and traditional mathematical modeling techniques.

\item[\dm] Highlighting Manopt's Enabling Role: This work underscores Manopt's crucial role in democratizing the use of \gls{mo}. By abstracting the complex geometric calculations and solver implementations, Manopt allows researchers and engineers to concentrate on the core task of defining the optimization problem within a familiar and powerful programming environment, thereby broadening the scope and applicability of \gls{mo} techniques in power systems.

\item[\dm] Foundation for Advanced Applications (Expansion Planning): We establish that the efficiency and accuracy of the \gls{mo}-based power flow solutions, particularly the \gls{bfs}-Manopt approach for \gls{epds}, provide a solid foundation for future extensions to more complex planning and operational problems, such as network reconfiguration using Master-Slave architectures.
\end{itemize}

The motivations and contributions of this paper are significant in advancing the understanding and application of \gls{mo} techniques in the context of modern power systems. The successful application of \gls{mo} to complex power flow problems in this paper, particularly enabled by the flexibility of tools like Manopt, clearly demonstrates its transformative potential for power systems optimization. To fully realize this potential and extend these benefits, continued research and development are essential.

\subsection{ORGANIZATION OF THE PAPER}

The remaining content is organized as follows. Section \ref{sec:II} provides an overview of manifold concepts, alternatives, and fundamental tools. Section \ref{sec:III} structures the methodology for formulating and solving optimization problems with manifolds. Section \ref{sec:IV} presents gradient descent-based algorithms adapted for manifold optimization. Section \ref{sec:V} demonstrates practical applications in power systems optimization, focusing on power flow analysis in EPDS and EPTS using Manopt. Section \ref{sec:BarriersMO} addresses integration challenges of manifold optimization techniques in power systems. Finally, Section \ref{sec:VII} draws conclusions and perspectives on manifold optimization for power systems applications, concluding with future research directions.

\section{MANIFOLD FUNDAMENTALS} \label{sec:II}

In this section, we start by highlighting different alternatives to the \gls{mo} technique. In the subsequent subsection, we explain the \gls{mo} framework. Finally, we catalog a list of manifolds found in many different real-world problems.

\subsection{ALTERNATIVES TO MANIFOLD OPTIMIZATION TECHNIQUE}

There are some alternatives for solving non-convex optimization problems, including \glspl{hem}, {\it convex relaxation} techniques, \gls{ml}-based algorithms, and {\it gradient-based} methods. {\bf HEMs}, such as \gls{ga}, \gls{pso}, \gls{gwo}, \gls{sa}, among others, are capable of performing a global search and are less likely to get trapped in local minima, therefore, being suitable to be applied to a wide range of problems without requiring gradient information, however, they present demerits of {\it a}) computationally intensive, often requiring many function evaluations, making them computationally expensive; {\it b}) lack of guarantees of convergence to the global optimum and can be slow to converge.

The merits of {\bf convex relaxation} techniques, such as \gls{sdp}, and \gls{ccp}, include {\it a}) rigorous framework for approximating non-convex problems into convex ones, and b) polynomial-time solvability. However, these techniques suffer scalability issues, rapidly becoming computationally infeasible for large-scale problems; moreover, the quality of the solution provided by \gls{sdp} and \gls{ccp} depends on how well the non-convex problem can be approximated by a convex one.

Besides,  {\bf \gls{ml}-based algorithms} can present impressive results since it can deal with large-scale problems, providing sub-optimal solutions. However, their feasibility in real-world scenarios is often limited. This limitation arises because many \gls{ml}-based algorithms require an offline training stage (particularly \glspl{nn} in supervised learning), utilizing data collected from real-world scenarios, however, it can be incompatible with the highly dynamic nature of power systems environments. The necessity for constant adaptation in these environments makes it challenging to rely on pre-trained models. Therefore, their practical application in power systems remains constrained by these real-world considerations.

Finally, {\bf gradient-based} methods, such as {\it gradient descent, Newton's method}, and {\it conjugate gradient}, represent a competitive alternative to the \gls{mo} approach for problems where gradient information is available, revealing strong local convergence properties in such scenarios. However, gradient-based methods can easily become trapped in local minima. As a substantial limitation, these methods require the \gls{of} to be differentiable, which cannot always be practical.

The key {\bf pros} and {\bf cons} of \gls{mo} over traditional optimization methods in power systems are summarized in Table \ref{table:comparison}, and include a) natural handling of non-convex constraints; b) geometric property exploitation; c) local linearization; d) versatility in handling constraints and symmetries; and e) high-dimensional data management: power systems often deal with high-dimensional data. Therefore, \gls{mo} methods can transform this data into a more manageable form, improving signal processing and resource allocation optimization.

\begin{table*}[ht!]
\caption{Comparison of \gls{mo} methods to its alternatives}\label{table:comparison}
\centering
\begin{tabularx}{\textwidth}{|>{\raggedright\arraybackslash}p{2.6cm}|>{\raggedright\arraybackslash}p{7.5cm}|>{\raggedright\arraybackslash}p{6.3cm}|}

\hline

\textbf{Method} & \textbf{Pros} & \textbf{Cons} \\

\hline\hline

\textbf{\gls{mo}} 

& 

{\it Handling Non-Convex Constraints}: \gls{mo} methods naturally handle non-convex constraints by treating the problem as an optimization over a smooth manifold.&  {\it Complexity in Implementation}: Implementing \gls{mo} methods can be complex due to the need for specialized knowledge in differential geometry and manifold theory.\\
\cline{2-3}
& 
{\it Exploiting Geometric Properties}: \gls{mo} methods leverage the geometric properties of the problem, handling constraints and symmetries (orthonormality, low rank, positivity, and invariance), allowing for more efficient solutions. 
&  

{\it Algorithmic Design Challenges}: The manifold constraint adds complexity to the algorithmic design and theoretical analysis.\\

\cline{2-3}

& {\it Local Linearization}: Manifolds are locally Euclidean, enabling linear optimization techniques in a more generalized form. & 
{\it Computational Overhead}: While efficient, \gls{mo} methods can still be computationally intensive, especially for high-dimensional problems.\\
\cline{2-3}

 & {\it High-Dimensional Data Management}: \gls{mo} methods can transform high-dimensional data into a more manageable form, improving optimization in tasks like signal processing and resource allocation. & \\
 
\hline\hline

\textbf{\glspl{hem}} 

& {\it Global Search Capability}: These methods perform a global search and are less likely to get trapped in local minima. 
&  

{\it Computationally Intensive}: They often require many function evaluations, making them computationally expensive. 
\\

\cline{2-3}

& {\it Flexibility}: They can be applied to various problems without requiring gradient information. 
&
{\it Lack of Guarantees}: No guarantee to converge to the global optimum or can converge slowly. \\

\hline\hline
\textbf{Convex Relaxation Techniques} & 
{\it Mathematical Rigor}: These methods provide a rigorous framework for approximating non-convex problems. & 
{\it Approximation Quality}: The quality of the solution depends on how well the non-convex problem can be approximated by a convex one. 
\\

\cline{2-3}

& {\it Polynomial-Time Solvability}: Convex problems can be solved efficiently using polynomial-time algorithms. & 
{\it Scalability Issues}: These methods can become computationally infeasible for large-scale problems. \\
\hline\hline

\textbf{\gls{ml}-based algorithms}
& \textit{Adaptability}: \gls{ml} methods can adapt to various scenarios and data patterns without requiring explicit modeling of the underlying physical processes.
& \textit{Training Data Requirement}: \gls{ml} methods require large amounts of high-quality training data, which may not always be available or easy to obtain.
\\
\cline{2-3}
& \textit{Data-Driven}: \gls{ml} methods leverage large datasets to learn and improve performance over time, making them suitable for environments where data is abundant.
& \textit{Computational Complexity}: Training \gls{ml} models, especially deep learning models, can be computationally intensive and time-consuming.
\\
\cline{2-3}
& \textit{Automation}: Once trained, \gls{ml} models can automate complex decision-making processes, reducing the need for manual intervention.
& \textit{Generalization}: \gls{ml} models may struggle to generalize well to unseen scenarios or out-of-distribution data, leading to suboptimal performance.
\\
\cline{2-3}
& \textit{Scalability}: \gls{ml} algorithms can handle high-dimensional data and scale well with the increasing complexity of power systems.
& \textit{Interpretability}: Particularly deep \gls{nn}, often act as black boxes, making it difficult to interpret/understand their decision-making processes.
\\

\hline\hline

\textbf{Gradient-Based Methods}
& {\it Efficiency}: These methods are efficient for problems where gradient information is available. 
& {\it Requirement of Smoothness}: These methods require the \gls{of} to be differentiable. \\
\cline{2-3}
& {\it Local Convergence}: They have strong local convergence properties.  & {\it Local Minima}: They can easily get trapped in local minima.\\
\hline
\end{tabularx}
\end{table*}

\subsection{MANIFOLD OPTIMIZATION FRAMEWORK}
In an optimization framework, we consider the search space $\mathcal{S}$ as the set containing all possible answers to our problem, and a cost function $f:\mathcal{S} \rightarrow \mathbb{R}$ which associates a cost $f(x)$ to each element $x$ of $\mathcal{S}$. The goal is to find $x \in \mathcal{S}$ such that $f(x)$ is minimized: 
\begin{small}
\begin{align}
\arg \min_{x\in \mathcal{S}} f(x).
\end{align}
\end{small}

We occasionally wish to denote the subset of $\mathcal{S}$ for which the minimal cost is achieved. We should bear in mind that this set might be empty. 

The Euclidean structure of $\mathbb{R}^n$ and the \gls{of} $f$'s smoothness are irrelevant to the optimization problem's definition. They are merely structures that we should use algorithmically to our advantage. Assuming linearity, the \gls{mo} approach requires smoothness as the key structure to exploit. 

\subsubsection{Optimization Over Smooth Surfaces}

Manifolds are a fundamental concept in mathematics, particularly in geometry and topology. Manifolds provide a generalization of shapes and spaces that locally resemble Euclidean space. In fact, \gls{mo} is a versatile framework for continuous optimization. It encompasses optimization over vectors and matrices and allows optimizing over curved spaces to handle constraints and symmetries such as orthonormality, low rank, positivity, and invariance under group actions \cite{manopt2014}.

Consider the set $\mathcal{M}$ as a smooth manifold, and the function $f$ is smooth on $\mathcal{M}$. Optimization over such surfaces can be understood as constrained because $x$ is not free to travel in $\mathbb{R}^n$ space, but is allowed only to stay on the surface. The favored alternative viewpoint, in this case, is to consider this as unconstrained optimization in a universe where the smooth surface is the only thing that exists. As a result, generalized Euclidean methods from unconstrained optimization can be applied to the larger class of optimization over smooth manifolds. We require a correct knowledge of gradient and Hessian on smooth manifolds to generalize techniques such as gradient descent and Newton's method. In the linear case, this requires including an inner product or an Euclidean structure. In a more general situation, it is advisable to exploit the property that smooth manifolds are locally linearizable around all points. The linearization at $x$ is the tangent space. Giving each tangent space its inner product\footnote{Varying smoothly with $x$ in a way to be determined precisely.} transforms the manifold into a \gls{rm}, upon which we construct what is known as a Riemannian structure \cite{Boumal_book_2023, manopt2014}.

\subsubsection{Operators on Riemannian Manifold}

\gls{rm}s are mathematical objects that generalize the notion of Euclidean space to more complex and curved geometries. These spaces are fundamental in various fields, including optimization, differential geometry, and theoretical physics \cite{Absil2007}. A Riemannian manifold is locally similar to an Euclidean space, but differs in that it is equipped with a Riemannian metric tensor. This tensor defines the distances and angles between points on the manifold by assigning a positive definite inner product to each tangent space. This inner product allows the measurement and interpretation of geometric properties such as length, angle, and curvature. Some key definitions and concepts in Riemannian geometry include:

\begin{itemize}
    \item[\dm]Riemannian Gradient ($\nabla_{\mathcal{M}} f(\boldsymbol{x})$): This is the generalization of the gradient from Euclidean space  $\nabla f$ to Riemannian manifolds. Specifically, the Riemannian gradient of a function $f$ on a manifold $\mathcal{M}$ is the projection of the Euclidean gradient onto the tangent space of the manifold at a given point.  
    \begin{small}
    \begin{align}
        \nabla_{\mathcal{M}} f(\boldsymbol{x}) = \textrm{Proj}_{\mathcal{T}_x \mathcal{M}} (\nabla f(\boldsymbol{x})),
    \end{align}
    \end{small}
    where $\textrm{Proj}_{\mathcal{T}_x \mathcal{M}}(\cdot)$ is the projection operator onto the tangent space $\mathcal{T}_x \mathcal{M}$.   It represents the direction of the steepest ascent of a function $f$ on the manifold $\mathcal{M}$.

\item[\dm]Retraction Operation ($\operatorname{Retr}_{\mathcal{M}}(\boldsymbol{x})$): The retraction operator $\operatorname{Retr}_{\mathcal{M}}(\boldsymbol{x})$ of a point on a manifold $\mathcal{M}$ is the projection of the given point $\boldsymbol{x}$ over the manifold $\mathcal{M}$. Retractions are used to ensure that optimization steps remain on the manifold.
\end{itemize}

The Riemannian gradient and retraction operation are essential for algorithms that optimize manifolds, as they ensure that the iterative steps respect the manifold's geometric structure. Moreover, we should bear in mind that each manifold has its own projection operator on the tangent space, as well as the retraction operator.

\subsubsection{Challenges in Manifold Optimization}

If additional constraints other than the manifold constraint are applied, one can add an indicator function of the feasible set of such additional constraints in the \gls{of}. Hence, the optimization problem covers a general formulation for \gls{mo}. Moreover, the manifold constraint is one of the main difficulties in algorithmic design and theoretical analysis.

One of the main challenges in \gls{mo} usually is the non-convexity of the manifold constraints. By utilizing the geometry of the manifold, a large class of constrained optimization problems can be viewed as unconstrained optimization problems on the manifold \cite{Hu2020}.

 \subsection{COLLECTION OF MANIFOLDS}

Optimization on manifolds is a versatile framework for continuous optimization. It encompasses optimization over vectors and matrices and adds the possibility to optimize over curved spaces to handle constraints and symmetries such as orthonormality, low rank, positivity, and invariance under group actions.

One of the most common manifolds is the {\bf \gls{ccm}}, in which all elements of the optimization variable must have a unit modulus. This naturally arises in voltage phasor optimization, where bus voltages are constrained to fixed magnitudes (e.g., in voltage-controlled (PV) buses or unity power factor inverters). Hence, the \gls{mo} framework is well-suited for power system problems involving phasor constraints, such as optimal power flow or distributed control of inverter-based resources.

Table \ref{tab:manifolds} summarizes the common real and complex types of manifolds, with particular emphasis on the \gls{ccm}, also known as the “complex one-manifold”.

A {\it complex manifold} is a manifold with a structure that locally resembles complex Euclidean space, {\it i.e.}, $\mathbb{C}^n$. This means a neighborhood is homeomorphic around every point to an open subset of $\mathbb{C}^n$. In particular, Table \ref{tab:ComplexCircle} shows the main features and applications of the Complex Circle $(\mathcal{S}^1)$ manifold.

\begin{table*}[!htbp]
\centering
 \caption{Common collection of manifolds \cite{Boumal_book_2023}}
\label{tab:manifolds}
\begin{tabular}{|p{1.7cm}|p{15cm}|} 
\hline
\bf Manifolds & \bf Feature \\
\hline\hline
\bf Euclidean Space $\mathbb{R}^n$ & $ \mathbb{R}^n $ is the most straightforward example of a manifold, where each point has a local neighborhood that looks exactly like $ \mathbb{R}^n$. Flat, infinite extent, commonly used in most basic analyses. \\
\hline
\bf Circle ($\mathcal{S}^1$) & $ \mathcal{S}^1 $ represents a one-dimensional manifold (1-manifold), which can be thought of as points equidistant from a center point in 2D space, like the perimeter of a circle; intrinsic periodicity (models cyclical phenomena). \\
\hline
\textbf{Sphere ($\mathcal{S}^n$)} & $\mathcal{S}^n$ generalizes the concept of a circle and sphere to “$n$” dimensions;  e.g., $ S^2 $ is the 2D surface of a 3D ball. Compact, without boundary, intrinsic higher-dimensional analogs. {\bf Use Cases:} Modeling surfaces like Earth's surface $(S^2)$. \\
\hline
\bf Torus $(\mathcal{T}^2)$ & The 2D torus is a surface shaped like a donut, which can be defined as $\mathcal{S}^1 \times \mathcal{S}^1$, the product of two circles or, generalizing, a product of $n$ circles, closed and compact. {\bf Use Cases:} Modeling periodic boundary conditions, complex cyclical phenomena \\
\hline
{\bf Projective Space} $\mathbb{RP}^n$ & Space of lines through the origin, compact, involves projective transformation. {\bf Use Cases:} Computer vision, robotics, projective geometry. \\
\hline
\textbf{Hyperbolic Space} $\mathbb{H}^n$ & Non-Euclidean, negatively curved. {\bf Use Cases:} Representing hierarchical tree structures, complex networks. \\
\hline
{\bf Hyperplanes} & These are generalizations of planes in higher dimensions. \\
\hline
\textbf{Lie Groups:} & Smooth manifold that is also a group, with applications in physics and engineering. \textit{Examples:} $SO(3)$, $SU(2)$. {\bf Use Cases:} Robotics, control theory, representation of symmetries. \\
\hline
\textbf{Grassmannian} ($G(k, n)$) & Space of all $ k $-dimensional subspaces of an n-dimensional vector space. {\bf Use Cases:} Signal processing, principal component analysis in higher dimensions. \\
\hline
\textbf{Stiefel Manifold} ($V(k, n)$): & Space of all orthonormal $ k $-frames in $ n $-space. {\bf Use Cases:} Multivariate statistics, optimization on orthonormal matrices. \\
\hline
\textbf{Kähler Manifold:} & A complex manifold with a Hermitian metric, deeply tied to complex and symplectic geometry. \textbf{Use Cases:} Theoretical physics, string theory. \\
\hline
\textbf{Calabi-Yau Manifold:} & A special type of Kähler manifold with a Ricci-flat metric. {\bf Use Cases:} String theory, particularly compactification methods. \\
\hline
\multicolumn{2}{|c|}{\textit{\textbf{Complex Manifolds}}} \\
\hline
\vspace{1.3mm}
\bf Complex Circle $(\mathcal{S}^1)$ & or {\it Complex 1-Manifold}: identified with the complex circle; defined as the set of all complex numbers of the unit norm, defined as: $\mathcal{S}^1 = \{ z \in \mathbb{C} \mid |z| = 1 \}$. Here, $|z|$ denotes the modulus of the complex number $z$. Manifolds modeled on complex numbers, allowing holomorphic coordinates. {\bf Use Cases:} Complex dynamics, algebraic geometry. \\
\hline
\end{tabular}
\end{table*}
\begin{table}[!ht]
\centering
\caption{Features and applications for the complex circle $(\mathcal{S}^1)$ Manifold}
\label{tab:ComplexCircle}
\begin{tabular}{|p{1.4cm}|p{6cm}|} 
\hline
\bf  Feature  & \bf  Description\\
\hline\hline
\textbf{1-Dim.} & While being embedded in $\mathbb{C}$ (which is like $\mathbb{R}^2$), the complex circle $\mathcal{S}^1$ is a 1-dimensional manifold.
\\
\hline
\textbf{Compactness} & It is a closed and bounded subset of $\mathbb{C}$.
\\
\hline
\textbf{Local Structure} & Locally, around any point on $\mathcal{S}^1$, it resembles the real line $\mathbb{R}$, meaning it can be mapped one-to-one onto an open interval of $\mathbb{R}$
\\
\hline
\bf Visualization &  One can visualize $\mathcal{S}^1$ as the unit circle in the complex plane, where each point on the circle is defined by a complex number $z$ with $|z| = 1$. This can be parameterized as $z = e^{i\theta}$ for $\theta \in [0, 2\pi)$, capturing its circular nature.
\\
\hline
\bf Complex Structure &  Looking at local coordinates as a complex manifold using complex logarithms and exponential. These give the local diffeomorphisms needed to open up parts of $\mathbb{C}$. The manifold structure is given by charts that map intervals around each point to the Euclidean space $ \mathbb{C}$. 
\\
\hline
\multicolumn{2}{|c|}{\textsc{Applications of the Complex Circle Manifold}}\\
\hline
\textbf{Topology}& Understanding the structure and properties of $\mathcal{S}^1$ is fundamental in algebraic topology, which contributes to studying fundamental groups and covering spaces. 
\\
\hline
\bf Complex Structure &
$\mathcal{S}^1$ is the natural domain for periodic functions and is central in studying Fourier analysis. 
\\
\hline
 \textbf{Physics}&  The complex circle appears in various physical theories, including quantum mechanics and wave mechanics, describing spaces of phases.
 \\
\hline
\end{tabular}
\end{table}

\section{METHODOLOGY FOR FORMULATING AND SOLVING OPTIMIZATION PROBLEMS WITH MANIFOLDS}\label{sec:III}

In the following, we illustrate the \gls{mo} methodology by detailing the steps involved in addressing a real-world problem in power systems. Specifically, we will optimize the voltage phase angles in a distributed grid using a \gls{mo} technique. This approach serves as an example, but the methodology can be applied to any problem involving non-convex constraints that can be represented as a manifold.

Power systems are designed to balance dynamically the generation and demand while maintaining stability; thus, optimizing voltage angles to enhance power flow efficiency is crucial, particularly in networks with high renewable energy penetration. The phase angles of bus voltages and the setpoints of distributed energy resources constitute high-dimensional and nonlinear spaces. \gls{mo} provides a structured framework for efficiently dealing with these complex and non-convex optimization problems.

\subsection*{STEP 1: PROBLEM FORMULATION}

\begin{itemize}

\item[\dm]\textbf{Define the Objective}: Clearly define the optimization problem in terms of an \gls{of}. This could involve optimizing key operational metrics such as power loss minimization, voltage stability margin, or economic dispatch efficiency. The \gls{of} is represented by $f(\boldsymbol{\theta},\boldsymbol{v})$, where $\boldsymbol{v}$ and $\boldsymbol{\theta}$ denote the complex bus voltages and phase angles, respectively, in a system with $N$ buses, with $[\boldsymbol{v}]_n = V_n e^{j\theta_n}$.

    \vspace{2mm}
    
    \item[\dm]\textbf{Constraints}: Identify the constraints of the problem, such as inferior and superior bounds on the injected power, represented as $\mathcal{P}_n^{\min} \leq P_n \leq \mathcal{P}_n^{\max}$ $\forall n \in \mathcal{N}_{\text{PV}}$, voltage magnitude bounds $V_n^{\min} \leq |V_n| \leq V_n^{\max}$ $\forall n \in \mathcal{N}_{\text{PV}}$, and power balance equations $\mathbf{S} = \mathbf{diag(V)YV^*}$, with $\mathbf{V}$ being the complex voltage and $\mathbf{Y}$ the adimittance matrix. For voltage-controlled buses, the following must hold:
    \begin{small}
    \begin{align}
    |V_n e^{j\theta_n}| = V_n^{\text{ref}} \quad \forall n \in \mathcal{N}_{\text{PV}},
    \end{align}
    \end{small}
    \noindent where $\mathcal{N}_{\text{PV}}$ denotes the set of PV buses and $V_n^{\text{ref}}$ is the specified voltage magnitude.
\end{itemize}

\subsection*{STEP 2: IDENTIFY THE MANIFOLD STRUCTURE}

\begin{itemize}
\item[\dm]\textbf{Manifold Description}: Establish the geometric structure of the constraint. For instance, voltage phasors in power systems lie on \gls{ccm}, which can be treated as a manifold\footnote{Common manifolds in power systems include the Euclidean manifold (for traditional power flow and OPF), the Grassmann manifold (for low-rank subspace methods in PMU data), and the SPD manifold (for covariance-based uncertainty analysis). For a complete list of manifolds, see Table \ref{tab:manifolds}.} and described as
    \begin{small}
    \begin{align}
        \mathcal{S}^1 = \{V_n e^{j\theta_n} \in \mathbb{C} \mid \theta_n \in [0, 2\pi), V_n = V_n^{\text{ref}}\}, \quad \forall n \in \mathcal{N}_{\text{PV}},
    \end{align}
    \end{small}
    \noindent where $\mathcal{N}_{\text{PV}}$ is the set of voltage-controlled (PV) buses in the system.
\end{itemize}

\subsection*{STEP 3: REFORMULATE THE OPTIMIZATION PROBLEM}

\begin{itemize}
    \item[\dm]\textbf{Manifold Representation}: Reformulate the original non-convex constraint problem in terms of manifold constraints. For example, if the optimization involves fixed-voltage-magnitude constraints, e.g., PV bus voltage angles $\theta_n \in [0, 2\pi)$, it should be represented in terms of the complex phasor form $V_ne^{j\theta_n}$. In power system optimization, the problem becomes:
    \begin{small}
    \begin{align}
        \max_{\boldsymbol{\theta} \in \mathcal{M}} f(\boldsymbol{\theta},\boldsymbol{V}),
    \end{align}
    \end{small}
    \noindent where $\boldsymbol{\theta}$ is the vector of voltage phase angles at PV buses, $\boldsymbol{V}$ represents the complex bus voltages, and $\mathcal{M}$ denotes the power flow manifold incorporating both voltage magnitude constraints and network physics.

    \vspace{2mm}
    
    \item[\dm]\textbf{Alternative Parameterization}: Use appropriate parameterizations to represent elements on the manifold in a computationally friendly way.
\end{itemize}

\subsection*{STEP 4: DEVELOP AN OPTIMIZATION ALGORITHM}

\begin{itemize}
    \item[\dm]\textbf{Initialization}: Start with an initial feasible point on the manifold. This might involve random initialization or a heuristic-based initialization.
    \begin{small}
    \begin{align}
        \boldsymbol{\theta}^{(0)} = [\theta_1^{(0)}, \theta_2^{(0)}, \ldots, \theta_N^{(0)}]^T.
    \end{align}
    \end{small}

    \vspace{2mm}
     
    \item[\dm]\textbf{Gradient Descent}: Use some {\it manifold-based optimization technique} (see subsection \ref{subsec:gradient} of section \ref{sec:V}) to iteratively update the phase shifts and move towards the local optimum. The steps should ensure the updated points lie on the manifold.
    \begin{small}
    \begin{align}
        \boldsymbol{\theta}^{(k+1)} = \boldsymbol{\theta}^{(k)} + \alpha \nabla_{\mathcal{M}} f(\boldsymbol{\theta}^{(k)}),
    \end{align}
    \end{small}
    \noindent where $\alpha$ is the step size and $\nabla_{\mathcal{M}} f(\boldsymbol{\theta}^{(k)})$ is the Riemannian gradient at the $k$-th iteration.

    \vspace{2mm}

    \item[\dm]{\bf Returning to the Manifold}: After updating the phase shifts, the updated point should be on the manifold surfaces, therefore, the retraction operator should be applied. Specifically, for the \gls{ccm} manifold, the retraction operator is given as
    \begin{small}
    \begin{align}
        \operatorname{Retr}_{\mathcal{S}^1}(\boldsymbol{\theta}^{(k+1)}) = \frac{\left[\boldsymbol{\theta}^{(k+1)}\right]_n}{\left|\left[\boldsymbol{\theta}^{(k+1)}\right]_n \right|}, \quad \forall n={1,2,\dots,N}. \label{eq:Retraction_CCM}
    \end{align}
    \end{small}
\end{itemize}

\subsection*{STEP 5: ITERATIVE OPTIMIZATION}

\begin{itemize}
    \item[\dm]\textbf{Iterative Process}: Iterate the optimization process until convergence. The stopping criterion could be based on the change in the \gls{of} value or the gradient norm at the $k$-th iteration till a small positive threshold $\epsilon$.
\end{itemize}
\begin{small}
\begin{align}
    \|\nabla_M f(\boldsymbol{\theta}^{(k)})\| < \epsilon.
\end{align}
\end{small}

\subsection*{STEP 6: VALIDATION AND PERFORMANCE EVALUATION}

\vspace{2mm}

\begin{itemize}
    \item[\dm]\textbf{Validation}: Validate the optimized voltage profiles by evaluating the system's power flow solutions. Compare the performance with conventional optimization methods (e.g., interior-point or Newton-Raphson) to demonstrate the efficiency and effectiveness of the \gls{mo}-based approach through metrics like:
    \begin{itemize}
        \item[\dm] Power loss reduction: $\Delta P_{\text{loss}} = P_{\text{loss}}^{\text{conv}} - P_{\text{loss}}^{\text{MO}}$
        \item[\dm] Voltage stability margin improvement: $\Delta VSM = VSM^{\text{MO}} - VSM^{\text{conv}}$
        \item[\dm] Computational time savings: $t_{\text{savings}} = t_{\text{conv}} - t_{\text{MO}}$
    \end{itemize}
\end{itemize}

\section{GRADIENT DESCENT-BASED ALGORITHMS} \label{sec:IV}
In the following, we exhibit a list of gradient descent-based methods, which can be generalized to the Riemannian space and utilized in \gls{mo} strategy.

\subsection{GRADIENT DESCENT-BASED ALGORITHMS} \label{subsec:gradient}

Manifold-based optimization techniques, such as  \gls{rn}, \gls{rcg}, \gls{rtr} and \gls{rgd} are powerful for optimizing functions constrained to manifolds. These methods are specifically designed to account for the non-Euclidean geometry of the optimization space. By leveraging the inherent geometric properties of these manifolds, such techniques enable efficient and effective optimization pathways. This is particularly beneficial in solving complex optimization problems found in power systems, where gradients can be computed directly on the manifold, allowing algorithms to navigate the curvature of the space effectively. This tailored approach helps ensure efficient convergence to locally optimal solutions within these constrained environments.

\vspace{2mm}

\noindent{\em \gls{rgd}}: The \gls{rgd} is an optimization technique that extends the classical gradient descent method to Riemannian manifolds. The optimization process is constrained to a curved space rather than a flat Euclidean space. The key idea is to iteratively move towards the local optimum while ensuring that each update remains on the manifold. The \gls{rgd} Algorithm steps are as follows:

\begin{enumerate}
    \item[\dm]\textbf{Initialization:} Start with an initial feasible point on the manifold.
    
    \item[\dm]\textbf{Gradient Computation:} Compute the Riemannian gradient, which is the projection of the Euclidean gradient onto the tangent space of the manifold.
    
    \item[\dm]\textbf{Update Rule:} Move in the direction of the Riemannian gradient by a step size, ensuring the updated point lies on the manifold. This often involves a retraction operation that maps the point back onto the manifold.
    
    \item[\dm]\textbf{Iteration:} Repeat the gradient computation and update steps until convergence.
\end{enumerate}

\begin{algorithm}
\label{alg:rgd}
\small
\caption{\gls{rgd} Algorithm}
\begin{algorithmic}[1]

\State \textbf{Input:} Initial point $\mathbf{x}_0$ on manifold $\mathcal{M}$, step size $\alpha$, max iterations $K$

\State \textbf{Output:} Optimized point $\mathbf{x}^*$ on the manifold
\algorithmiccomment{\it Input and output specifications}

\State $k \gets 0$ \algorithmiccomment{ Initialization of the iteration counter}

\While{$k < K$ and not converged} 
    \State Compute Riemannian gradient $\operatorname{grad} f(\mathbf{x}_k)$ \algorithmiccomment{\it Project Euclidean gradient onto tangent space}

    \State Update: $\mathbf{x}_{k+1} \gets \mathcal{R}_{\mathbf{x}_k}(-\alpha \operatorname{grad} f(\mathbf{x}_k))$ 
    \algorithmiccomment{\it Utilize retraction operator to map point back to manifold} 
    
    \State $k \gets k + 1$
\EndWhile

\State \textbf{return} $\mathbf{x}_k$ \algorithmiccomment{\it  Return of the final optimized point}

\end{algorithmic}
\end{algorithm}

\vspace{2mm}

\noindent{\em \gls{rn}:} \gls{rn} method is an extension of the classical Newton method to Riemannian manifolds. It uses second-order information (Hessian) to achieve faster convergence when compared to the gradient descent method.

\vspace{2mm}

\noindent{\em \gls{rcg}:} \gls{rcg} is an adaptation of the conjugate gradient method to Riemannian manifolds. It combines the efficiency of the conjugated gradient method with the geometric constraints of manifolds.

\vspace{2mm}

\noindent{\em \gls{rtr}:} \gls{rtr} methods extend trust-region methods to the Riemannian manifolds. These methods iteratively solve a local approximation of the optimization problem within a “trust region” around the current point.

 Utilizing \gls{mcs} or other simulation techniques is useful to evaluate the algorithm's performance under different scenarios. Based on simulation results, we can also refine the algorithm for better performance by adjusting parameters and re-evaluating. Refinement based on simulations is a critical step in the optimization process, especially when dealing with complex systems like those involving manifolds. In this context, the refinement process based on simulations might involve specific steps such as:

\begin{itemize}
    \item[\dm] \textit{Gradient and Hessian Adjustments:} Fine-tuning the computation of Riemannian gradients and Hessians to ensure they accurately capture the manifold's geometry.
    \item[\dm] \textit{Retraction Operations:} Modifying the retraction operations to better map updated points back onto the manifold.
    \item[\dm] \textit{Step Size Adaptation:} Adjusting the step size dynamically based on the manifold's curvature to ensure efficient convergence.
    \item[\dm] \textit{Constraint Handling:} Refining how additional constraints are incorporated into the optimization problem, possibly by adjusting the indicator functions or penalty terms.
\end{itemize}

The goal is to iteratively improve the algorithm or model based on empirical evidence from simulations, leading to better performance in the target application. How refinement is typically done based on simulations is summarized in Algorithm \ref{alg:refinement}. 

Some open-source libraries support manifold structures, such as  \href{https://www.manopt.org/}{\texttt{Manopt tools}}, which implements a bunch of manifold collections available on \texttt{Manopt} (Matlab) or \texttt{Pymanopt} (Python), which can be useful for obtaining solutions of many diverse problems \cite{manopt2014}, 
\cite{Boumal_book_2023}.

\begin{algorithm}
\caption{Refinement Based on Simulations}
\label{alg:refinement}
\small
\begin{algorithmic}[1]
    \State \textbf{Initial Simulation:} 
    \State \quad Run initial simulations to evaluate the performance of the current algorithm or system configuration.
    \State \quad Collect performance metrics such as \gls{se}, \gls{ee}, error rates, throughput, etc.
    
    \State \textbf{Performance Analysis:}
    \State \quad Examine the simulation results to identify strengths and weaknesses.
    \State \quad Determine which aspects are underperforming or causing issues.
    
    \State \textbf{Parameter Adjustment:}
    \State \quad Modify parameters of the algorithm or system based on the analysis.
    \State \quad Make necessary changes to the algorithm itself if required.
    
    \State \textbf{Re-simulation:}
    \State \quad Perform new simulations with the adjusted parameters or modified algorithm.
    \State \quad Compare the new simulation results with the previous ones.
    
    \State \textbf{Iteration:}
    \State \quad Repeat the process of analysis, adjustment, and re-simulation multiple times.
    \State \quad Continue iterating until performance metrics converge to satisfactory levels.
    
    \State \textbf{Validation:}
    \State \quad Use cross-validation techniques to ensure generalization to different scenarios.
    \State \quad Test the refined system under various conditions for robustness.
    
    \State \textbf{Final Tuning:}
    \State \quad Perform fine-tuning of parameters to achieve the best possible performance.
    \State \quad Use advanced optimization techniques if necessary.
\end{algorithmic}
\end{algorithm}


\section{MANIFOLDS IN POWER SYSTEMS OPTIMIZATION} \label{sec:V}

In the sequel, we present detailed applications of the \gls{mo} framework to power system optimization problems in \gls{epds} and \gls{epts}. Specifically, we demonstrate how different \gls{mo} techniques can be effectively employed for solving power flow. Through the application of suggested \gls{mo} algorithms and methodologies, we verify the significant performance improvements achievable when employing modern \gls{mo} approaches in power systems analysis.

\subsection{\bf POWER FLOW IN ELECTRIC POWER DISTRIBUTION SYSTEMS}

In \gls{epds}, the power flow assumes specific characteristics that must be considered. To model the power flow problem for \gls{epds}, the following conditions are taken into account:

\begin{itemize}
    \item[\dm] The system is represented by a single-line equivalent diagram per phase;
    \item[\dm] Loads are modeled as constant demanded powers;
    \item[\dm] The substation is the only power source;
    \item[\dm] Active and reactive power losses in the lines are concentrated at the initial bus of the transmission line;
    \item[\dm] The capacitive reactance of the lines is not considered.
\end{itemize}

\begin{figure}
\begin{center}
\includegraphics[width=8cm]{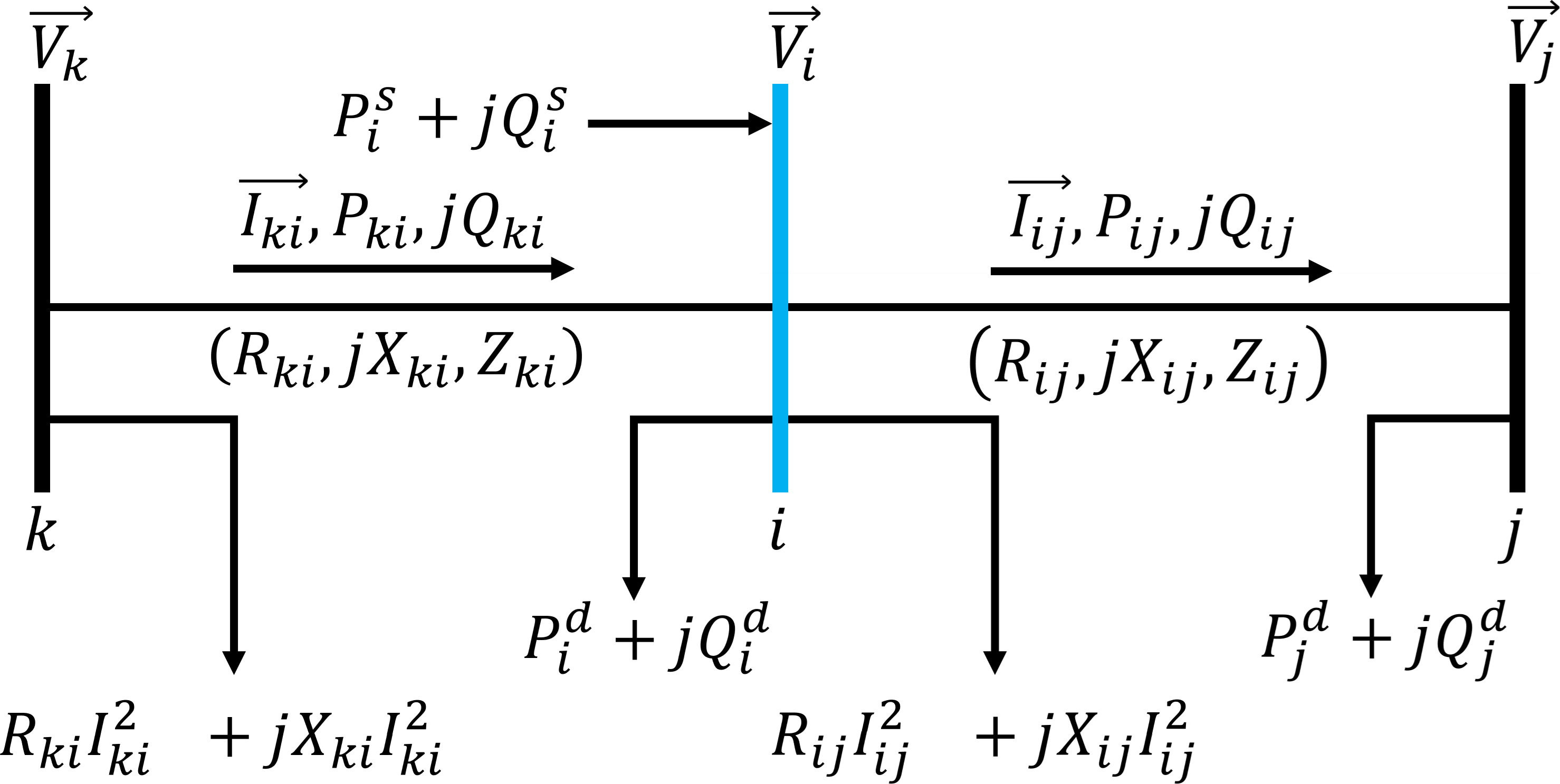}
\caption{Power flow diagram for \gls{epds}.}
\label{FP:fig1}
\end{center}
\end{figure}

Thus, Figure~\ref{FP:fig1} presents the power flow diagram for the \gls{epds}, which is the basis for mathematically modeling power flow.

Therefore, we aim to minimize the \gls{of} in Equation~\eqref{FP:eq1}, whose single component is the power losses, subject to Equations~\eqref{FP:eq2}–\eqref{FP:eq7}.

\begin{small}
\begin{align}
    &\underset{I_{ij},P_{ij},Q_{ij}~\forall ij \in \Omega_l;~V_i,P^s_i,Q^s_i~\forall i \in \Omega_b}{\text{minimize}}~f = \sum_{\forall ij \in \Omega_l} R_{ij} I^2_{ij}, \label{FP:eq1} \\
    &\text{subject to:} \nonumber \\
    &\sum_{\forall ki \in \Omega_l} P_{ki}-\sum_{\forall ij \in \Omega_l} (P_{ij} + R_{ij} \cdot I^2_{ij}) + P^s_i = P^d_i;~ \forall_i \in \Omega_b \label{FP:eq2} \\
    &\sum_{\forall ki \in \Omega_l} Q_{ki}-\sum_{\forall ij \in \Omega_l} (Q_{ij} + X_{ij} \cdot I^2_{ij}) + Q^s_i = Q^d_i;~ \forall_i \in \Omega_b \label{FP:eq3} \\
    &V^2_i - 2(R_{ij} P_{ij} + X_{ij} Q_{ij}) - Z^2_{ij} I^2_{ij} - V^2_j = 0;~ \forall ij \in \Omega_l \label{FP:eq4} \\
    &I^2_{ij} V^2_j = P^2_{ij} + Q^2_{ij};~ \forall ij \in \Omega_l \label{FP:eq5} \\
    &\underline{I}_{ij} \leq I_{ij} \leq \overline{I}_{ij};~ \forall ij \in \Omega_l \label{FP:eq6} \\
    &\underline{V}_i \leq V_i \leq \overline{V}_i;~ \forall i \in \Omega_b, \label{FP:eq7}
\end{align}
\end{small}

The variables $I_{ij}$, $P_{ij}$, and $Q_{ij}$ represent the current, active power and reactive power in the lines, $V_i$ represents the voltage at the buses and $P^s_i$ and $Q^s_i$ represent the active and reactive power from the substation. The parameters $R_{ij}$, $X_{ij}$ and $Z_{ij}$ represent the resistances, inductive reactances, and impedances of the lines, $\overline{I}_{ij}$ and $\underline{I}_{ij}$ represent the upper and lower limits for the current, $\overline{V}_i$ and $\underline{V}_i$ represent the upper and lower limits for the voltage and $P^d_i$ and $Q^d_i$ represent the demanded active and reactive power.

The conventional active and reactive power balance equations are shown in Equation~\eqref{FP:eq2}-\eqref{FP:eq3}, respectively, where we can see that the power demanded at a bus $i$ $P^d_i$ must be equal to the power that $i$ receives from $k$ $P_{ki}$ and the powers distributed from $i$ to $j$ $P_{ij}$ plus the losses $R_{ij} \cdot I^2_{ij}$. If the bus is a substation, $P_{ki}$ is zero and $P^s_i$ is nonzero.

To define the voltage drops, considering the circuit represented in Figure~\ref{FP:fig1} and the fact that $\vec{V}_i - \vec{V}_j = \vec{Z}_{ij} \vec{I}_{ij} = (R_{ij} + jX_{ij}) \vec{I}_{ij};~ \forall ij \in \Omega_l$, we have Equation~\eqref{FP:eq4}.

Also, the apparent power constraint in each branch of the network is given by Equation~\eqref{FP:eq5}, the current limit of the lines is given by Equation~\eqref{FP:eq6}, and the voltage limit at the buses is provided by Equation~\eqref{FP:eq7}.

It is possible to note that the \gls{of} of Equation~\eqref{FP:eq1} is separable, and by definition, the sum of convex functions is a convex function. Therefore, as Equation~\eqref{FP:eq8} suggests, given that the Hessian of each term in the summation is positive definite, the \gls{of} is convex.

\begin{small}
\begin{align}
    &H\left(R_{ij} I^2_{ij}\right) = 2R_{ij} \succ 0 \label{FP:eq8}
\end{align}
\end{small}

However, it may be problematic to evaluate the geometry of the constraint set and verify whether it would be convex. An alternative would be to use the CVX modeling system \cite{CVX} and check if \gls{dcp} rules are violated. It is possible to note that $I_{ij}$, $V_i$, and $V_j$ are squared in all constraints, except in those that define limits for them, which makes the expressions in Equations~\eqref{FP:eq1}-\eqref{FP:eq4} nonlinear. Additionally, Equation~\eqref{FP:eq5} has an invalid quadratic form according to DCP rules, since the product of $I^2_{ij}$ and $V^2_j$ would not be a square even if the variables were not squared. Thus, the constraint set is not convex due to these quadratic equalities.

In the literature, the power flow model and the strategies derived from it are often linearized, involving discretizations of nonlinear terms. However, this approach implies that the accuracy of the solution inevitably depends on parameters tuned based on the solution of the original nonlinear and non-convex model, which is not ideal, in addition to introducing uncertainties.

Convexifying the power flow problem into an \gls{socp} model has positive implications for the planning and expansion of electrical systems. Strategies of combinatorial nature that are derived directly from the power flow, such as network reconfiguration, can inherit the convexity of the base model. These strategies rarely introduce new nonlinearities and are generally based on binary decision variables, so the relaxation of the Big-M method is usually sufficient for cases where there is a product of continuous and binary variables.

However, both the linearized and convex models face the same problem regarding the choice of solver, as solvers that support integer variables are often commercial: AMPL \cite{AMPLIDE}, a modeling editor widely used in power system modeling, includes SCIP, a free solver for \gls{mip} and \gls{minlp} problems that can also solve \gls{socp} problems, but the AMPL environment itself is not free to use; CVX \cite{CVX} only provides Gurobi, GLPK, and MOSEK as solvers with support for integer constraints, which are commercial solvers; Pyomo \cite{Pyomo} is slower and memory-inefficient, and although it has free solvers like CBC or HiGHS that can solve MILP problems, the uncertainties introduced by linear models persist and the time spent adjusting parameters remains a problem, besides, it is not possible to solve an \gls{socp} problem with these solvers.

Nonlinear solvers traditionally target smooth problems with continuous variables, however, some solvers accept integer variables, and some aim for global optimality. The few global nonlinear solvers that can handle integer variables, such as BARON, LINDO, RAPOSa, and COUENNE available in AMPL, are commercial and expensive solvers, and the rest of the local solvers ensure that their solution satisfies the first-order KKT conditions, which is usually not enough.

Therefore, solvers that solve \gls{mip} and \gls{minlp} are rarely non-commercial, and furthermore, the tendency is that these non-commercial solvers have a higher execution time, as is the case, for example, with SDPT3 in CVX compared to Gurobi. It is safe to say that the mathematical modeling of power flow has as one of its main objectives to enable the use of solvers commonly used in the literature, such as CPLEX, Knitro, and others, to solve combinatorial problems, whether nonlinear, linear, or convex, and if this were not possible, the mathematical modeling of power flow might not be so interesting.

It is expected that commercial solvers promote very low execution time, and it may not be interesting to look for some exception to try to solve combinatorial problems, as the tendency is for the complexity of these models to increase as new research is published, execution time and even accuracy may end up falling short. An alternative would be the use of a Master-Slave configuration, in which the Master sends a suggestion of the main variable of the problem, usually binary or integer, so that the Slave returns the power flow calculation.

\begin{figure}
\begin{center}
\includegraphics[width=6cm]{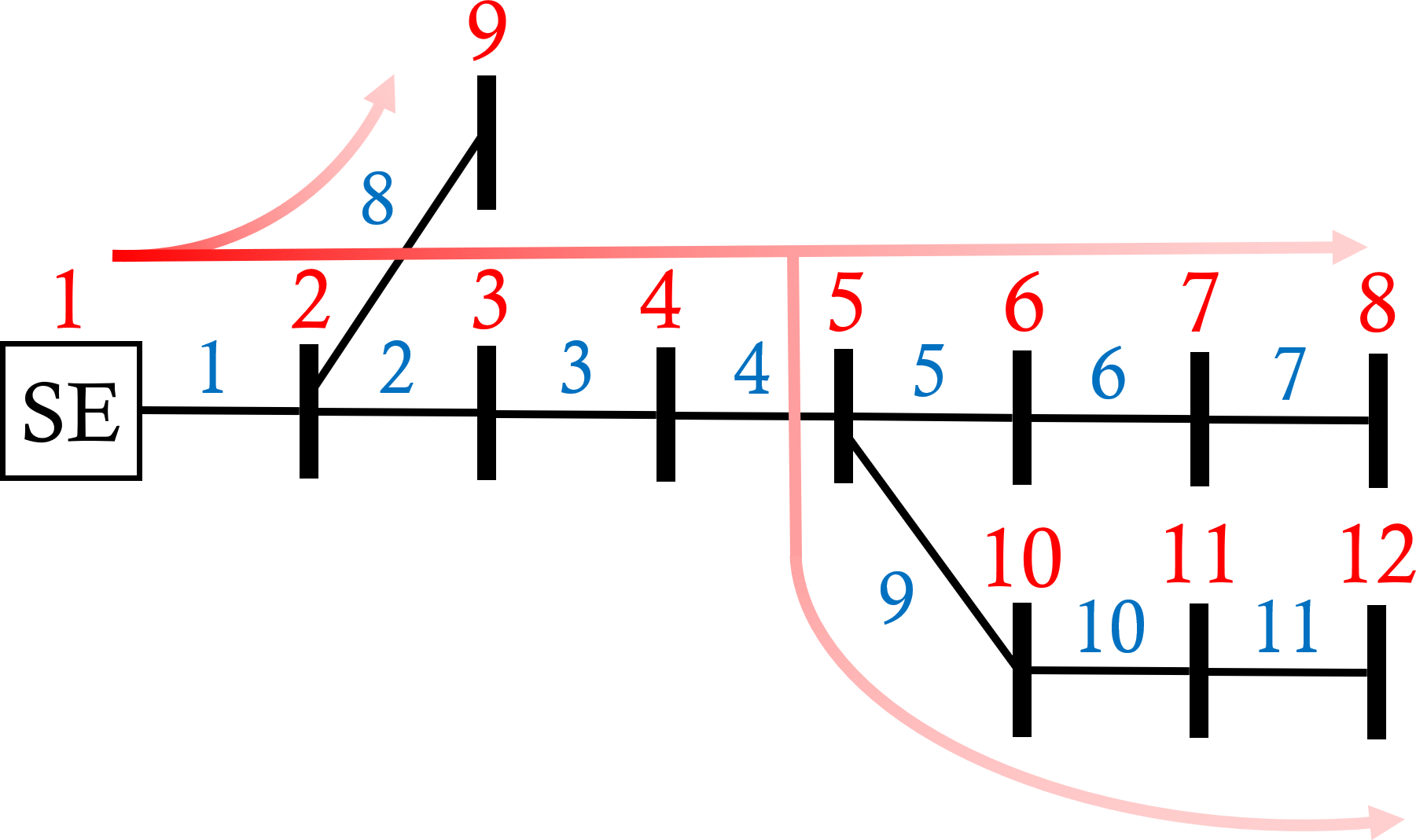}
\caption{Expected power flow in a 12-bus \gls{epds}.}
\label{FP:fluxo}
\end{center}
\end{figure}

Thus, in the context of manifold and in terms of computational time, the model presented earlier may have problems regarding the choice of the initial point, considering that, as Figure~\ref{FP:fluxo} suggests, the many variables of the problem must inevitably be computed and may ultimately depend on how power flows in the system. Since shared buses across branches lead to cumulative power and current calculations, it becomes necessary to enumerate the branches and perform computations resembling a \gls{bfs} approach. Specifically, the calculation of power relies on currents to account for losses, while the determination of currents depends on voltages through the apparent power equation, and voltage calculations require currents to estimate voltage drops. This interdependence necessitates an iterative process of successive calculations and recalculations to converge upon a solution.

The \gls{bfs} algorithm is a widely used method for power flow analysis in \gls{epds}, particularly in radial and weakly meshed systems. Studies have demonstrated its computational efficiency and robustness compared to traditional methods like Newton-Raphson and Zbus. For instance, \cite{Zimmerman1996} reported that \gls{bfs}, especially when combined with network reduction and fast decoupling, executes power flow analyses more than three times faster than NR and Zbus methods. BFS and its variants, such as \gls{fbs}, Adaptive \gls{bfs}, and \gls{qbbfs}, have been extensively evaluated in diverse network conditions. \cite{Muruganantham2016} highlighted \gls{fbs}'s superior performance in high resistance-to-reactance scenarios, where Newton-Raphson methods often fail to converge. Similarly, \cite{Mahmoud2016} proposed \gls{qbbfs} for multi-phase systems, noting its robustness and lower iteration counts compared to standard \gls{bfs}. In unbalanced low-voltage networks, \cite{GarcesRuiz2022} observed that \gls{bfs}, though requiring more iterations than Newton's method, consistently converged under challenging conditions.

It is reasonable to imagine, then, that a \gls{bfs} algorithm could be used instead of the presented modeling to avoid all the problems mentioned, with the latter serving as a comparison. However, in the context of manifolds, it is necessary to derive a cost function from the \gls{bfs}.

\subsubsection{\bf MANIFOLD SELECTION AND RIEMANNIAN GRADIENT PROJECTION} 

In Manopt manifold-based optimization platform \cite{manopt2014}, the \textit{problem} structure is a fundamental part of defining an optimization problem on a manifold. It contains several fields describing the manifold structure, cost function, Euclidean and Riemannian gradient, Euclidean and Riemannian hessian, etc., but the initial concern would be to define a cost function and its gradient. Having in mind, first of all, an optimization model, since Manopt does not support constraints, we must transform the problem with constraints into an equivalent problem without constraints, and penalty methods are very popular for this task.

Methods for solving a constrained optimization problem in $n$ variables and $m$ constraints can be divided roughly into four categories:
\begin{itemize}
    \item[\dm] Primal methods: work in $n - m$ space,
    \item[\dm] Penalty methods: work in $n$ space,
    \item[\dm] Dual and cutting plane methods: work in $m$ space,
    \item[\dm] Lagrangian methods: work in $n + m$ space.
\end{itemize}

\noindent{\bf Penalty Method}: The classical penalty method was first developed by \cite{fiacco1968nonlinear}. The goal of the penalty approach is to convert the problem into an equivalent unconstrained problem or into a problem with simple bounded constraints (as unconstrained methods): a term is added to the \gls{of} to penalize any violation of the constraints, thus it generates a sequence of infeasible points whose limit is an optimal solution to the original problem. In this way, it is possible to solve nonlinear programming problems having equality and inequality constraints \cite{Bazaraa2006}. Consider the following problem:
\begin{small}
\begin{align}
\begin{array}{ll}
\underset{x}{\text{minimize} }& f(x), \\
\text{s.t.} & x \in S,
\end{array}
\end{align}
\end{small}
where $f(x)$ is a continuous function on $E^n$ and $S$ is a constraint set in $E^n$. In most applications, $S$ is defined implicitly by a number of functional constraints. The idea of a penalty function method is to replace the original problem by an unconstrained problem of the form:
\begin{small}
\begin{align}
\begin{array}{ll}
\underset{x}{\text{minimize}} & \theta(x) = f(x) + cP(x),
\end{array}
\end{align}
\end{small}

\noindent where $\theta(x)$ is the augmented \gls{of}, $c$ is a parameter that determines the relative importance of the constraints, and $P(x)$ is a function to measure violations of the constraints satisfying:
\begin{small}
\begin{align}
\begin{array}{ll}
i) & P(x) \text{ is continuous}, \\
ii) & P(x) \geq 0, \, \forall x \in E^n, \\
ii) & P(x) = 0, \, \text{if and only if } x \in S.
\end{array}
\end{align}
\end{small}

A suitable penalty function must incur a positive penalty for infeasible points and no penalty for feasible points. If the constraints are of the form:
\begin{small}
\begin{align}
\begin{array}{ll}
g_i(x) \leq 0, & i = 1, \ldots, m \quad \text{and} \quad h_i(x) = 0, \quad i = 1, \ldots, l,
\end{array}
\end{align}
\end{small}

\noindent a suitable penalty function $\alpha(x)$ is defined by
\begin{small}
\begin{align}
\alpha(x) = \sum_{i=1}^{m} \phi [g_i(x)] + \sum_{i=1}^{l} \psi [h_i(x)],
\end{align}
\end{small}

\noindent where $\phi(\cdot)$ and $\psi(\cdot)$ are continuous functions satisfying:
\begin{small}
\begin{align}
\begin{array}{ll}
\phi(y) = 0 & \text{if} \quad y \leq 0 \quad \text{and} \quad \phi(y) > 0 \quad \text{if} \quad y > 0 \\
\psi(y) = 0 & \text{if} \quad y = 0 \quad \text{and} \quad \psi(y) > 0 \quad \text{if} \quad y \neq 0 \\
\end{array}
\end{align}
\end{small}

Typically, $\phi$ and $\psi$ are of the forms:
\begin{small}
\begin{align}
\phi(y) = [\max\{0,~y\}]^p \quad \text{and} \quad \psi(y) = |y|^p, \label{phipsi}
\end{align}
\end{small}

\noindent where $p$ is a positive integer. Therefore, The penalty function $\alpha$ is usually of the form:
\begin{small}
\begin{align}
\alpha(x) = \sum_{i=1}^{m} [\max\{0,~g_i(x)\}]^p + \sum_{i=1}^{l} |h_i(x)|^p, \label{FP:alpha}
\end{align}
\end{small}

The augmented \gls{of} then becomes:
\begin{small}
\begin{align}
\theta(x) = f(x) + \mu \alpha(x).
\end{align}
\end{small}

Regarding the \gls{bfs}, the objective is generally to minimize the error related to supplying all demands as $\max \left( \left| S_D - V_{calc} \cdot I_n^* \right| \right)$. However, to reduce the number of variables that Manopt would need to handle, it is possible to minimize the difference between calculated and current voltages, thus making the only variables the voltage magnitudes and angles. The starting point does not need to be computed, and is simply a vector of nominal voltages with the imaginary part equal to zero. Also, having the goal of reaching zero in the cost function, this refers to the equality constraints. In this way, Algorithm~\ref{alg:bfs_power_flow} proposes a cost function that, using $p = 2$, takes the form:
\begin{small}
\begin{align}
\text{cost} &=  \sum_{i=1}^{l} |h_i(x)|^2. \label{fcost}
\end{align}
\end{small}

\begin{algorithm}
\small
\caption{BFS Power Flow Cost Function}
\label{alg:bfs_power_flow}
\begin{algorithmic}[1]
\State \textbf{Input:} 
\State \quad $x$ - State vector (voltage magnitudes and angles)
\State \quad $V_{nom}$ - Nominal voltage at slack bus
\State \quad $S_D$ - Apparent power demand at each bus
\State \quad $N_i$, $N_j$ - Sender and Receiver buses
\State \quad $Z$ - Branch impedances
\State \quad $N_b$ - Number of buses

\State \textbf{Output:} 
\State \quad $cost$ - Sum of squared residuals

\State $V(1) \gets V_{nom}$ \Comment{Slack bus voltage (considered to be bus 1 for convenience)}

\For{$k \gets 1$ \textbf{to} $N_b-1$}
    \State $V(k+1) \gets x(2k-1) + 1j \cdot x(2k)$ \Comment{Construct complex voltage}
\EndFor

\State \Comment{Calculate injected currents}
\For{$k \gets 1$ \textbf{to} $N_b$}
    \State $I_n(k) \gets (S_D(k) / V(k))^*$
\EndFor

\State \Comment{Backward sweep: calculate branch currents}
\State $I_n^{aux} \gets I_n$
\For{$i \gets \text{length}(N_i)$ \textbf{downto} $1$}
    \State $I_L(i) \gets I_n^{aux}(N_j(i))$
    \State $I_n^{aux}(N_i(i)) \gets I_n^{aux}(N_i(i)) + I_L(i)$
\EndFor

\State \Comment{Forward sweep: calculate voltages}
\State $V_{calc}(1) \gets V_{nom}$
\For{$i \gets 1$ \textbf{to} $\text{length}(N_i)$}
    \State $V_{calc}(N_j(i)) \gets V_{calc}(N_i(i)) - Z(i) \cdot I_L(i)$
\EndFor

\State \Comment{Difference between calculated and current voltages}
\State $residuals \gets V_{calc} - V$

\State \Comment{Cost: sum of squared residuals}
\State $cost \gets \text{sum}(\text{abs}(residuals)^2)$
\end{algorithmic}
\end{algorithm}

\vspace{2mm}\noindent\noindent{\bf Euclidean {\it vs.} Riemannian Gradient}: To derive the Riemannian gradient, we first need to choose the manifold. Since there are no additional restrictions beyond what is shown in Algorithm~\ref{alg:bfs_power_flow}, the Euclidean space may be sufficient. In this case, the optimization occurs in the space $\mathbb{R}^{2(n-1)}$, where $n$ represents the dimension of the optimization space, which corresponds to the number of buses in the \gls{epds}. The state of the power system is given by:
\begin{small}
\begin{align}
x &= [\operatorname{Re}\left(\boldsymbol{V}_{T_b=0}\right),~\operatorname{Im}\left(\boldsymbol{V}_{T_b=0}\right)] \in \mathbb{R}^{2(n-1)} \nonumber \\
&= [{\left(|V| \cos(\theta)\right)}_{T_b=0},~{\left(|V| \sin(\theta)\right)}_{T_b=0}] \in \mathbb{R}^{2(n-1)},
\end{align}
\end{small}

\noindent where $T_b$ is the type of bus (0 for load buses and 1 for the substation) and $\boldsymbol{V}_{T_b=0}$ is the state vector containing the complex form of voltage magnitudes and angles considering only load buses, since the substation voltage equals the system's nominal voltage. In accordance with \cite{FP:ref2}, \cite{FP:ref3} and \cite{FP:ref4}, we may hence define the set:
\begin{small}
\begin{align}
\mathcal{M} = \{ x \in \mathbb{R}^{2(n-1)} \mid h(x) = 0 \}, \label{manifold1}
\end{align}
\end{small}

\noindent where $h(x)$ encodes the model.

Thus, the projection of the gradient from the Euclidean space to the Riemannian space is trivial. To properly account for the absolute value of the complex-valued functions $h_i(x)$, we need to rewrite the cost function and its derivative. Since each $h_i(x)$ has real and imaginary components, we can express it as:
\begin{small}
\begin{align}
h_i(x) = h^\text{Re}_i(x) + j h^\text{Im}_i(x),
\end{align}
\end{small}

\noindent where $h^\text{Re}_i(x)$ and $h^\text{Im}_i(x)$ are the real and imaginary parts, respectively, and $j$ is the imaginary unit. The squared magnitude of $h_i(x)$ is:
\begin{small}
\begin{align}
|h_i(x)|^2 = h^\text{Re}_i(x)^2 + h^\text{Im}_i(x)^2.
\end{align}
\end{small}

Thus, the cost function becomes:
\begin{small}
\begin{align}
\text{cost} &= \sum_{i=1}^{l} |h_i(x)|^2 = \sum_{i=1}^{l} \left( h^\text{Re}_i(x)^2 + h^\text{Im}_i(x)^2 \right). \label{fcost_complex}
\end{align}
\end{small}

Next, to compute the gradient of the cost function with respect to $x$, we differentiate the squared magnitudes:
\begin{small}
\begin{align}
\frac{\partial}{\partial x} |h_i(x)|^2 &= \frac{\partial}{\partial x} \left( h^\text{Re}_i(x)^2 + h^\text{Im}_i(x)^2 \right)\nonumber \\
&= 2 h^\text{Re}_i(x) \nabla h^\text{Re}_i(x) + 2 h^\text{Im}_i(x) \nabla h^\text{Im}_i(x).
\end{align}
\end{small}

This can be written compactly as:
\begin{small}
\begin{align}
\frac{\partial}{\partial x} |h_i(x)|^2 = 2 \operatorname{Re}\left( h_i(x)^* \nabla h_i(x) \right),
\end{align}
\end{small}

\noindent where $h_i(x)^*$ is the complex conjugate of $h_i(x)$, and $\nabla h_i(x)$ is the gradient of $h_i(x)$ (treated as a complex-valued function).

The gradient of the total cost is then:
\begin{small}
\begin{align}
\nabla \text{cost} &= \sum_{i=1}^{l} \frac{\partial}{\partial x} |h_i(x)|^2 = 2 \sum_{i=1}^{l} \operatorname{Re}\left( h_i(x)^* \nabla h_i(x) \right). \label{FP:grad_complex}
\end{align}
\end{small}

\subsubsection{\bf RESULTS WITH THE MANOPT TOOLBOX}

\begin{figure}
\begin{center}
\includegraphics[width=7cm]{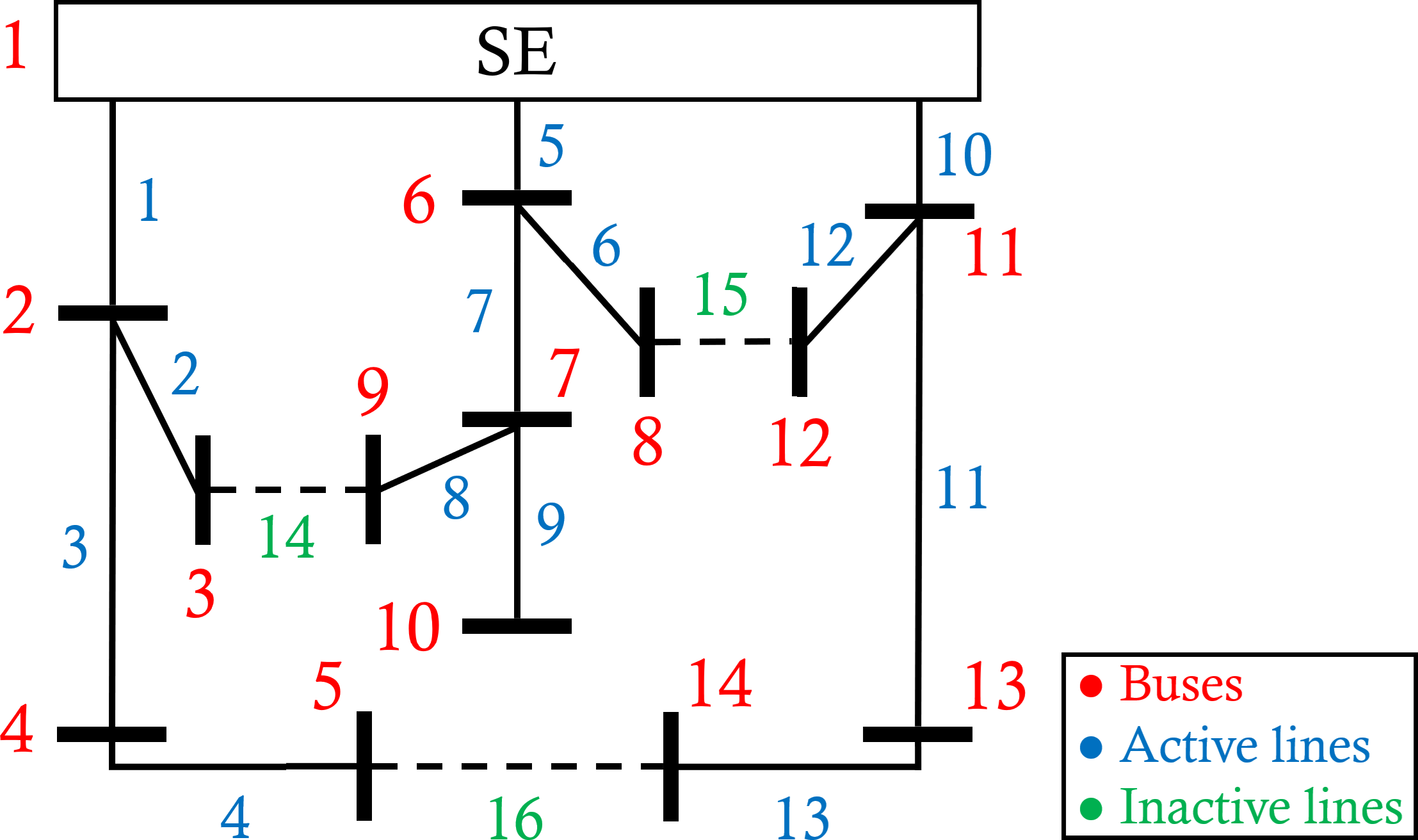}
\caption{14-bus \gls{epds}.}
\label{FP:fig14bussys}
\end{center}
\end{figure}

Using the 14-bus \gls{epds} shown in Figure~\ref{FP:fig14bussys}, which includes 1 installed substation, 13 load buses, a nominal voltage of 23 kV and a demand of (28700 + j5900) kVA with little reactive power demanded thanks to the action of capacitors, in accordance with the data presented by \cite{github}, the power flow can be solved using Manopt.

Therefore, we also want to compare the solution with the solution of other solvers for linear, nonlinear and convex models. The model presented at the beginning of this chapter is non-linear, the convex model can be written based on the relaxations presented by \cite{convex1} and \cite{convex2}, and the linear model can be written based on \cite{FP:ref6}. Table~\ref{FP:tab14} shows the comparison between power flow models for the 14-bus \gls{epds}.

\begin{table*}
\centering
\resizebox{1.5\columnwidth}{!}{%
\begin{tabular}{ccccc}
\hline \hline
    \textbf{Model} & \textbf{Options} & \textbf{Power losses (Kw)} & \textbf{Vmin (p.u.)} & \textbf{Time (s)} \\ \hline \hline
    \cite{19} & - & 511.4 & 0.9693 & - \\ \hline
    \cite{91} & - & 511.4 & 0.9693 & - \\ \hline
    \cite{67} & - & 511 & - & - \\ \hline
    \cite{76} & - & 511.436 & 0.96927 & - \\ \hline
    \cite{31} & - & 514 & - & - \\ \hline
    \cite{68} & - & 511.44 & 0.969 & 0.27 \\ \hline
    \cite{systems} & - & 511.44 & - & - \\ \hline
    \cite{7bus12bus} & - & 511.4 & 0.9693 & - \\ \hline
    \cite{FP:ref1} & - & 511.43 & 0.9692 & - \\ \hline
    \cite{92} & - & 511.43 & - & - \\ \hline
    Nonlinear (Knitro solver and & - & 511.443520 & 0.9693 & 0.131000 \\
    AMPL interface) & & & & \\ \hline
    Linear (CPLEX solver and & S = 6 & 511.464940 & 0.9693 & 0.066000 \\
    AMPL interface) & Y = 120 & & & \\ \hline
    Convex (CPLEX solver and & - & 511.443830 & 0.9693 & 0.045000 \\
    AMPL interface) & - & & & \\ \hline
    Linear (CPLEX solver and & S = 6 & 511.464940 & 0.9693 & 0.06 \\
    Pyomo interface) & Y = 120 & & & \\ \hline
    Convex (Gurobi solver and & - & 511.443520 & 0.9693 & 0.10 \\
    CVX interface) & & & & \\ \hline
    Convex (SDPT3 solver and & - & 511.443678 & 0.9693 & 0.50 \\
    CVX interface) & & & & \\ \hline
    Nonlinear BFS (steepestdescent & - & 511.443520 & 0.9693 & 0.146952 \\
    solver and Manopt interface) & & & & \\ \hline \hline
\end{tabular}}
\caption{Comparison between power flow models for the 14-bus \gls{epds}.}
\label{FP:tab14}
\end{table*}

\begin{figure}
\begin{center}
\includegraphics[trim=56pt 0pt 65pt 25pt,clip,width=\linewidth]{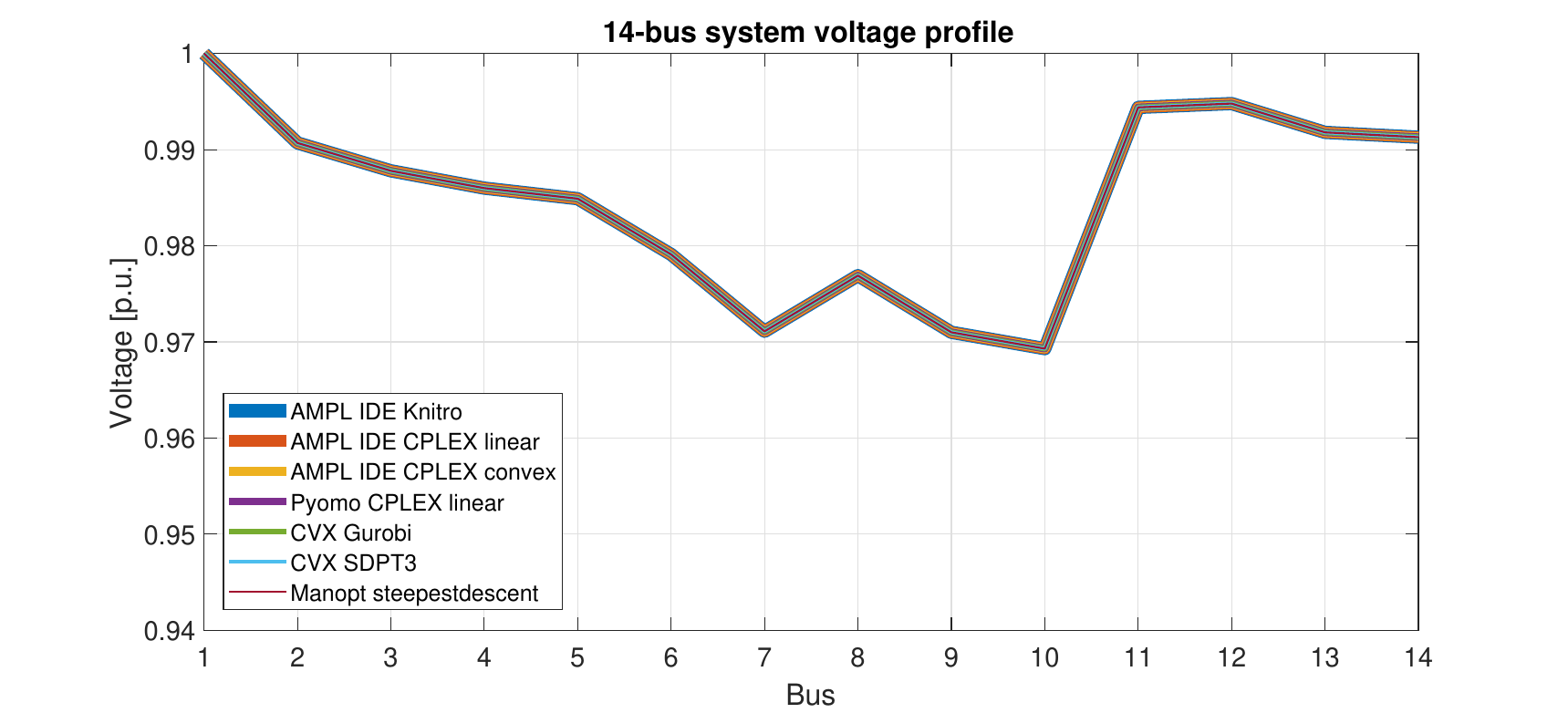}
\caption{Voltage profile for the 14-bus \gls{epds} power flow models.}
\label{FP:fig14busvolt}
\end{center}
\end{figure}

\begin{figure*}
\begin{center}
\includegraphics[trim=86pt 0pt 85pt 25pt,clip,width=.85\linewidth]{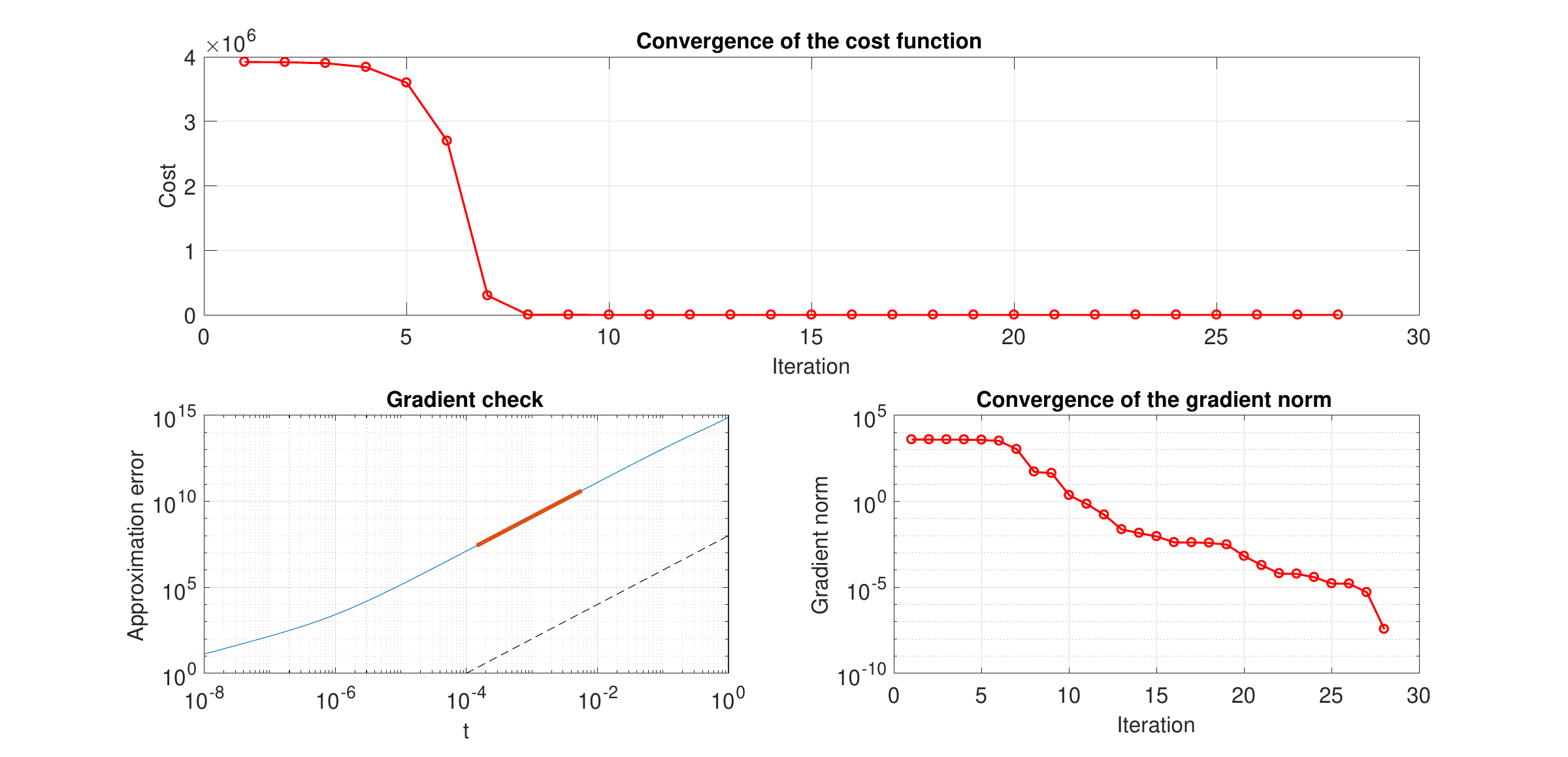}
\caption{Manopt metrics for the 14-bus \gls{epds} power flow \gls{bfs} model.}
\label{FP:manopt14}
\end{center}
\end{figure*}

Furthermore, Figure~\ref{FP:fig14busvolt} shows the voltage profiles. Since all the profiles overlap, it is safe to assume that the models were well-formulated and well-conditioned. Also, Figure~\ref{FP:manopt14} presents the Manopt metrics, showing the accuracy of the gradient, as well as the convergence of the cost function after 28 iterations and the gradient norm.
\includegraphics[,width=\linewidth]{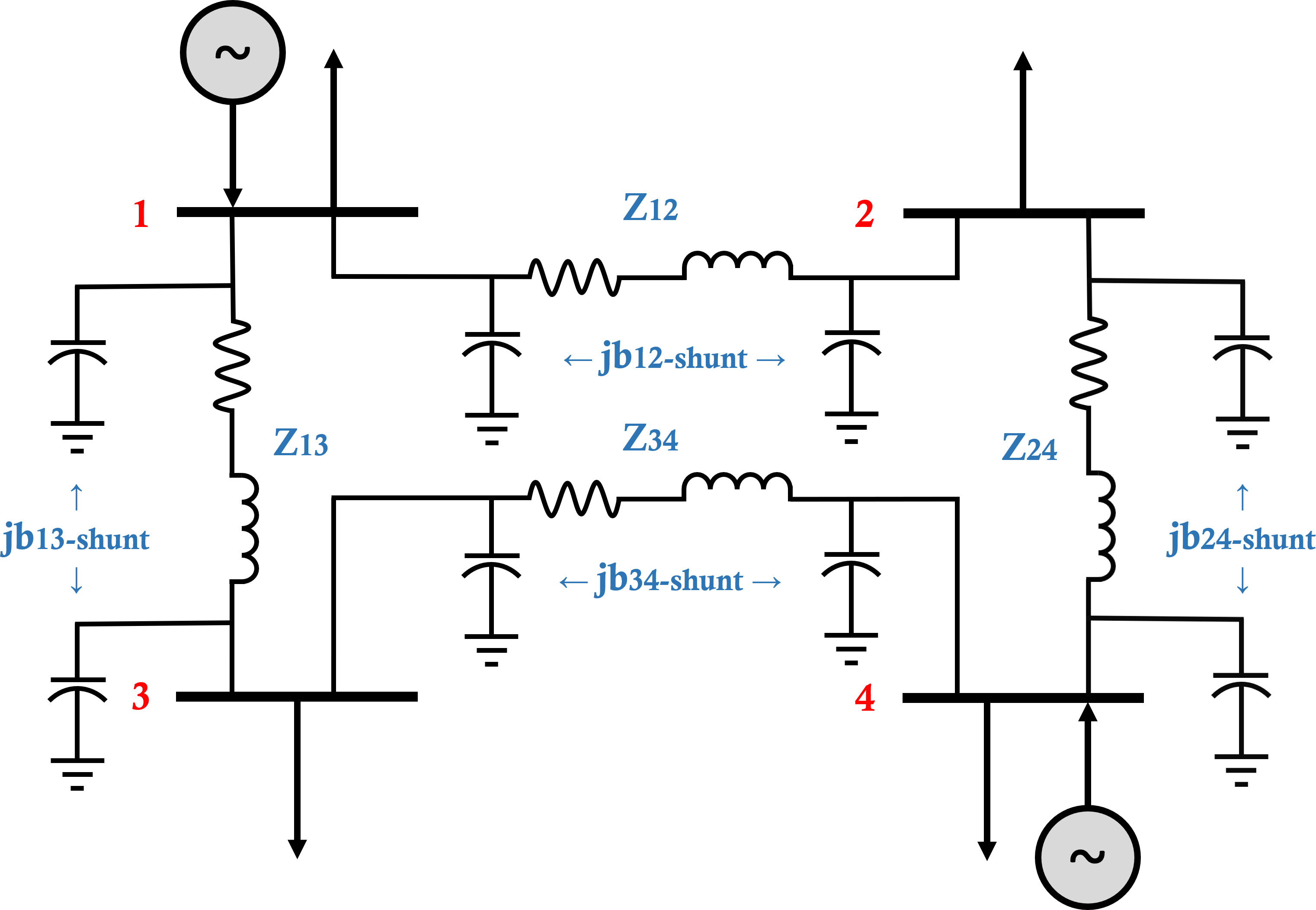}

Thus, the gradient provided is accurate. Notably, while Knitro, CPLEX and Gurobi are well-established commercial solvers and achieved very good results, Manopt stands out as a free and open-source alternative that matches its performance in terms of execution time and accuracy. The only models that exactly reached the value of the Knitro nonlinear model were Gurobi with the convex model and Manopt, and this is particularly impressive given the complexity of the nonlinear problem. Manopt achieved better results than CVX with SDPT3, which is also an interesting alternative in terms of freely distributed solvers, but not in terms of computational time, and it was also better than CPLEX in terms of accuracy.

\begin{figure}
\begin{center}
\includegraphics[width=8.4cm]{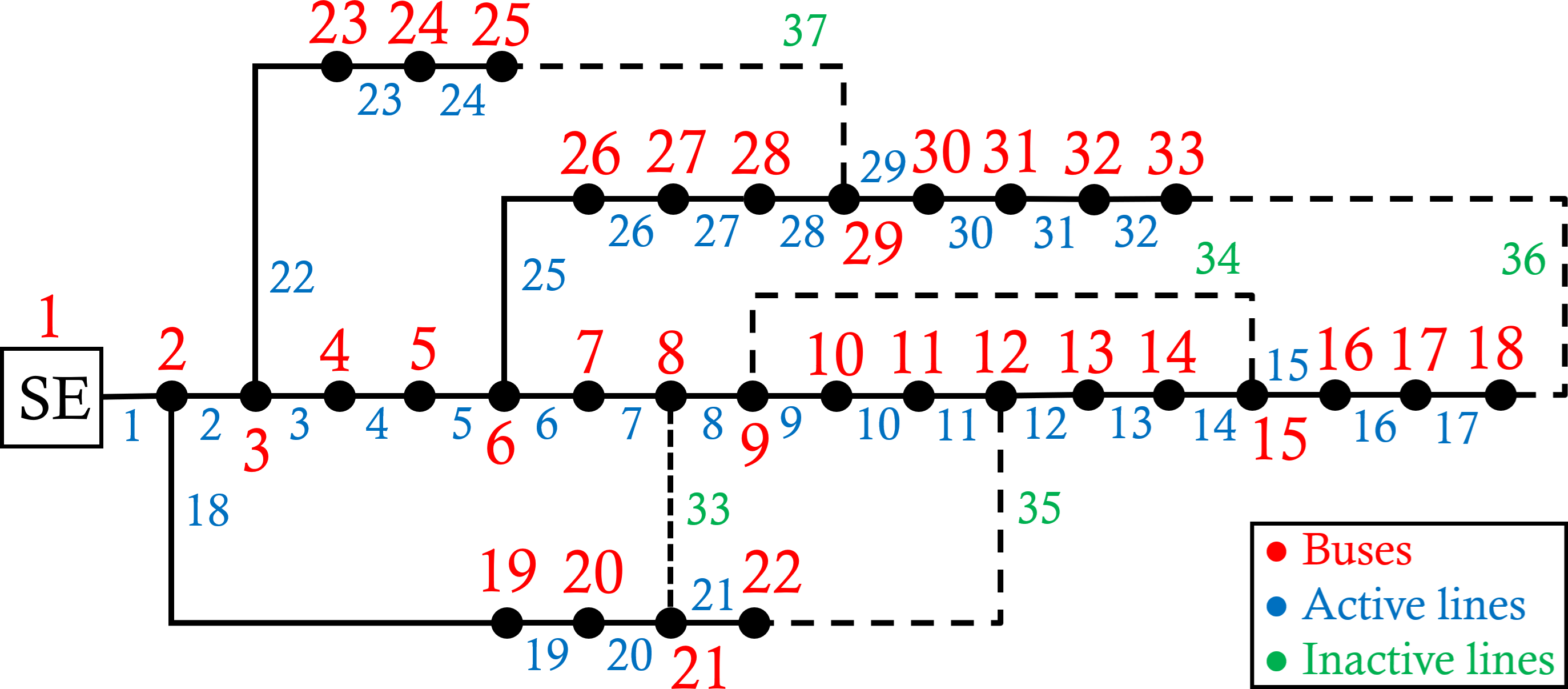}
\caption{33-bus \gls{epds}.}
\label{FP:fig33bussys}
\end{center}
\end{figure}

Subsequently, in order to verify another test case, using the 33-bus \gls{epds} shown in Figure~\ref{FP:fig33bussys}, which includes 1 installed substation, 32 load buses, a nominal voltage of 12.66 kV and a demand of (3715 + j2300) kVA, in accordance with the data presented by \cite{github}, the power flow can also be solved using Manopt.

Table~\ref{FP:tab33} shows the comparison between power flow models for the 33-bus \gls{epds}.

\begin{table*}
\centering
\resizebox{1.4\columnwidth}{!}{%
\begin{tabular}{ccccc}
\hline \hline
    \textbf{Model} & \textbf{Options} & \textbf{Power losses (Kw)} & \textbf{Vmin (p.u.)} & \textbf{Time (s)} \\ \hline \hline
    \cite{14} & - & 204.14 & 0.91 & - \\ \hline
    \cite{30} & - & 202.68 & - & - \\ \hline
    \cite{67} & - & 202.68 & - & - \\ \hline
    \cite{31} & - & 202.7 & - & - \\ \hline
    \cite{26} & - & 202.7 & - & - \\ \hline
    \cite{16} & - & 210.97 & 0.9037 & - \\ \hline
    \cite{systems} & - & 202.68 & - & - \\ \hline
    \cite{FP:ref1} & - & 202.56 & 0.9130 & - \\ \hline
    \cite{FP:ref6} & - & 202.56 & 0.9130 & - \\ \hline
    \cite{92} & - & 202.67 & - & - \\ \hline
    \cite{FP:ref5} & - & 202.67 & 0.9131 & - \\ \hline
    Nonlinear (Knitro solver and & - & 202.677126 & 0.9131 & 0.140000 \\
    AMPL interface) & & & & \\ \hline
    Linear (CPLEX solver and & S = 8 & 202.677029 & 0.9131 & 0.157000 \\
    AMPL interface) & Y = 78 & & & \\ \hline
    Convex (CPLEX solver and & - & 202.677416 & 0.9131 & 0.055000 \\
    AMPL interface) & - & & & \\ \hline
    Linear (CPLEX solver and & S = 8 & 202.677029 & 0.9131 & 0.17 \\
    Pyomo interface) & Y = 78 & & & \\ \hline
    Convex (Gurobi solver and & - & 202.677126 & 0.9131 & 0.21 \\
    CVX interface) & & & & \\ \hline
    Convex (SDPT3 solver and & - & 202.677143 & 0.9131 & 0.65 \\
    CVX interface) & & & & \\ \hline
    Nonlinear BFS (steepestdescent & - & 202.677126 & 0.9131 & 0.179256 \\
    solver and Manopt interface) & & & & \\ \hline \hline
\end{tabular}}
\caption{Comparison between power flow models for the 33-bus \gls{epds}.}
\label{FP:tab33}
\end{table*}

\begin{figure}
\begin{center}
\includegraphics[trim=56pt 0pt 65pt 25pt,clip,width=\linewidth]{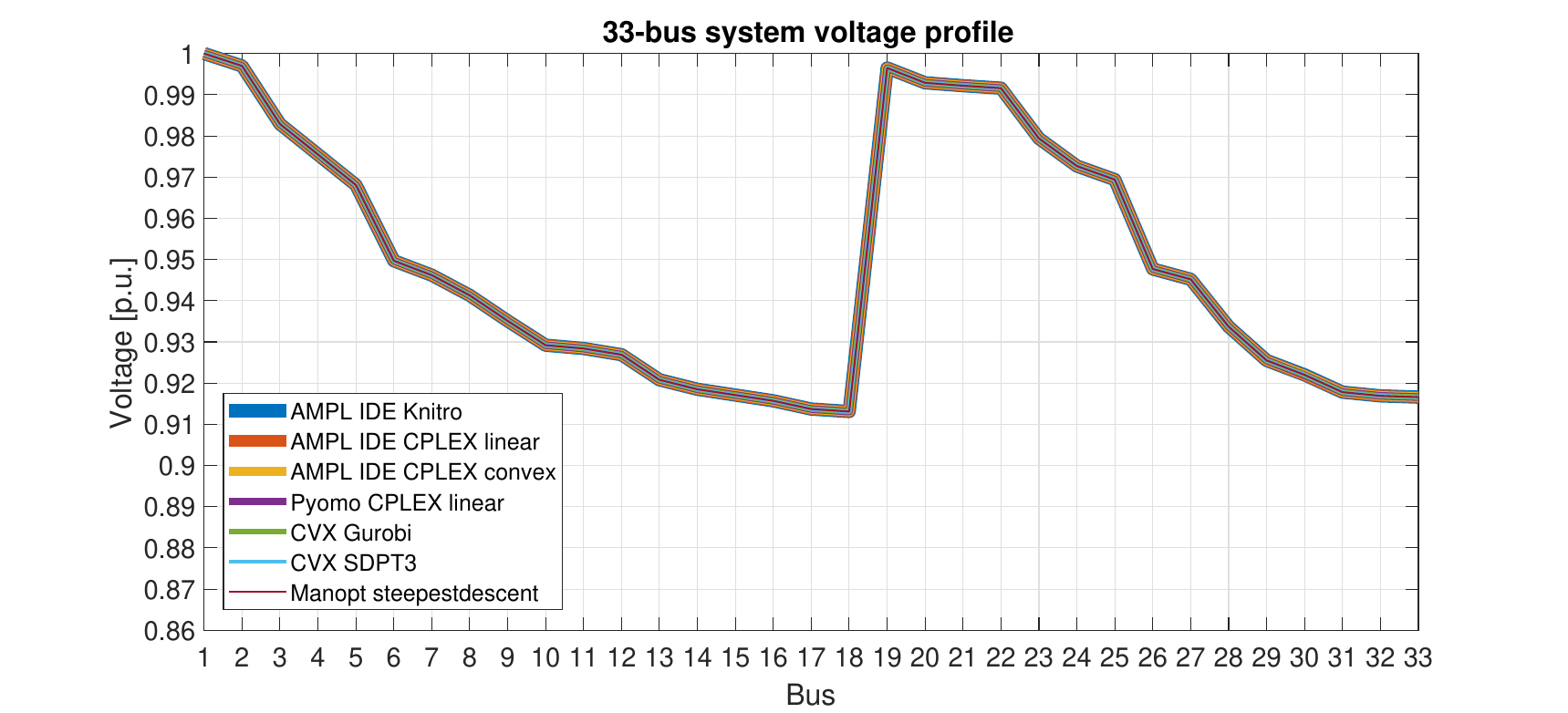}
\caption{Voltage profile for the 33-bus \gls{epds} power flow models.}
\label{FP:fig33busvolt}
\end{center}
\end{figure}

\begin{figure*}
\begin{center}
\includegraphics[trim=86pt 0pt 85pt 25pt,clip,width=.82\linewidth]{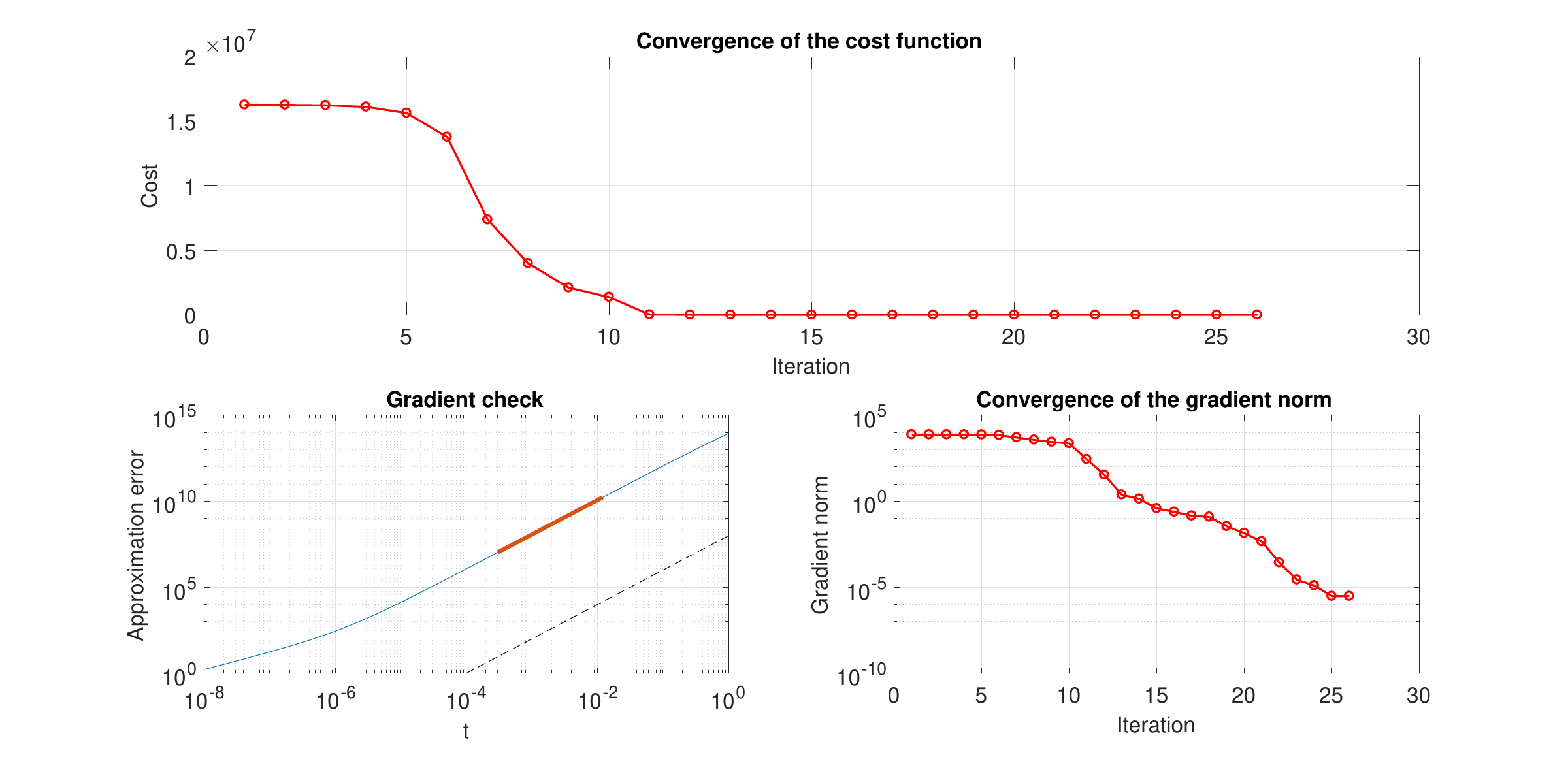}
\caption{Manopt metrics for the 33-bus \gls{epds} power flow \gls{bfs} model.}
\label{FP:manopt33}
\end{center}
\end{figure*}

Furthermore, Figure~\ref{FP:fig33busvolt} shows the voltage profiles. Since all the profiles overlap, it is again safe to assume that the models were well-formulated and well-conditioned. Also, Figure~\ref{FP:manopt33} presents the same Manopt metrics, showing the accuracy of the gradient, as well as the convergence of the cost function after 26 iterations and the gradient norm.

Thus, the gradient provided is accurate. Once again, the only models that exactly matched the value of the Knitro nonlinear model were Gurobi with the convex model and Manopt. In this case, Manopt stood out for its speed, being faster than Gurobi. In any case, the computational time between these results tends to be of the same order and very competitive. Once more, Manopt achieved better results than CVX with SDPT3 and CPLEX in terms of accuracy.

\begin{figure}
\begin{center}
\includegraphics[width=8cm]{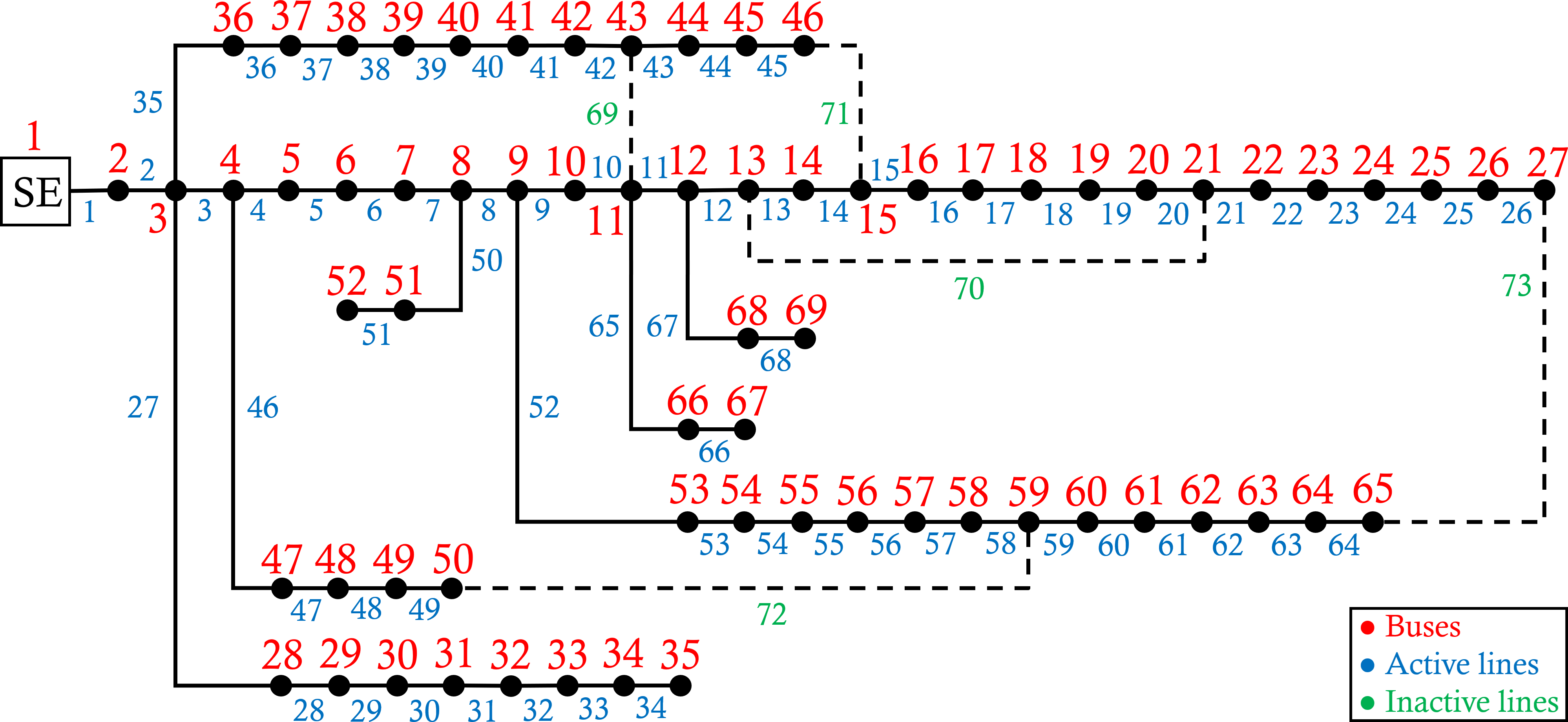}
\caption{69-bus \gls{epds}.}
\label{FP:fig69bussys}
\end{center}
\end{figure}

Finally, in order to verify one last test case, using the 69-bus \gls{epds} shown in Figure~\ref{FP:fig69bussys}, which includes 1 installed substation, 68 load buses, a nominal voltage of 12.66 kV and a demand of (3802.19 + j2694.60) kVA, in accordance with the data presented by \cite{github}, the power flow can also be solved using Manopt.

Table~\ref{FP:tab69} shows the comparison between power flow models for the 69-bus \gls{epds}.

\begin{table*}
\centering
\resizebox{1.55\columnwidth}{!}{%
\begin{tabular}{ccccc}
\hline \hline
    \textbf{Model} & \textbf{Options} & \textbf{Power losses (Kw)} & \textbf{Vmin (p.u.)} & \textbf{Time (s)} \\ \hline \hline
    \cite{30} & - & 225.00 & - & - \\ \hline
    \cite{68} & - & 224.93 & 0.969 & 20.2 \\ \hline
    \cite{13} & - & 224.7893 & 0.9092 & 0.359579 \\ \hline
    \cite{16} & - & 225 & 0.9091 & - \\ \hline
    \cite{systems} & - & 225 & - & - \\ \hline
    \cite{7bus12bus} & - & 225 & - & - \\ \hline
    \cite{FP:ref1} & - & 224.57 & 0.9091 & - \\ \hline
    \cite{FP:ref6} & - & 224.57 & 0.9091 & - \\ \hline
    \cite{FP:ref5} & - & 224.99 & 0.9092 & - \\ \hline
    \cite{18} & - & 224.9931 & 0.9092 & - \\ \hline
    Nonlinear (Knitro solver and & - & 224.993112 & 0.9092 & 0.122000 \\
    AMPL interface) & & & & \\ \hline
    Linear (CPLEX solver and & S = 7 & 224.978545 & 0.9092 & 0.237000 \\
    AMPL interface) & Y = 62 & & & \\ \hline
    Convex (CPLEX solver and & - & 224.993410 & 0.9092 & 0.075000 \\
    AMPL interface) & - & & & \\ \hline
    Linear (CPLEX solver and & S = 7 & 224.978545 & 0.9092 & 0.34 \\
    Pyomo interface) & Y = 62 & & & \\ \hline
    Convex (Gurobi solver and & - & 224.993113 & 0.9092 & 0.16 \\
    CVX interface) & & & & \\ \hline
    Convex (SDPT3 solver and & - & 224.993200 & 0.9092 & 0.69 \\
    CVX interface) & & & & \\ \hline
    Nonlinear BFS (steepestdescent & - & 224.993112 & 0.9092 & 0.266440 \\
    solver and Manopt interface) & & & & \\ \hline \hline
\end{tabular}}
\caption{Comparison between power flow models for the 69-bus \gls{epds}.}
\label{FP:tab69}
\end{table*}

\begin{figure}
\begin{center}
\includegraphics[trim=56pt 0pt 65pt 25pt,clip,width=\linewidth]{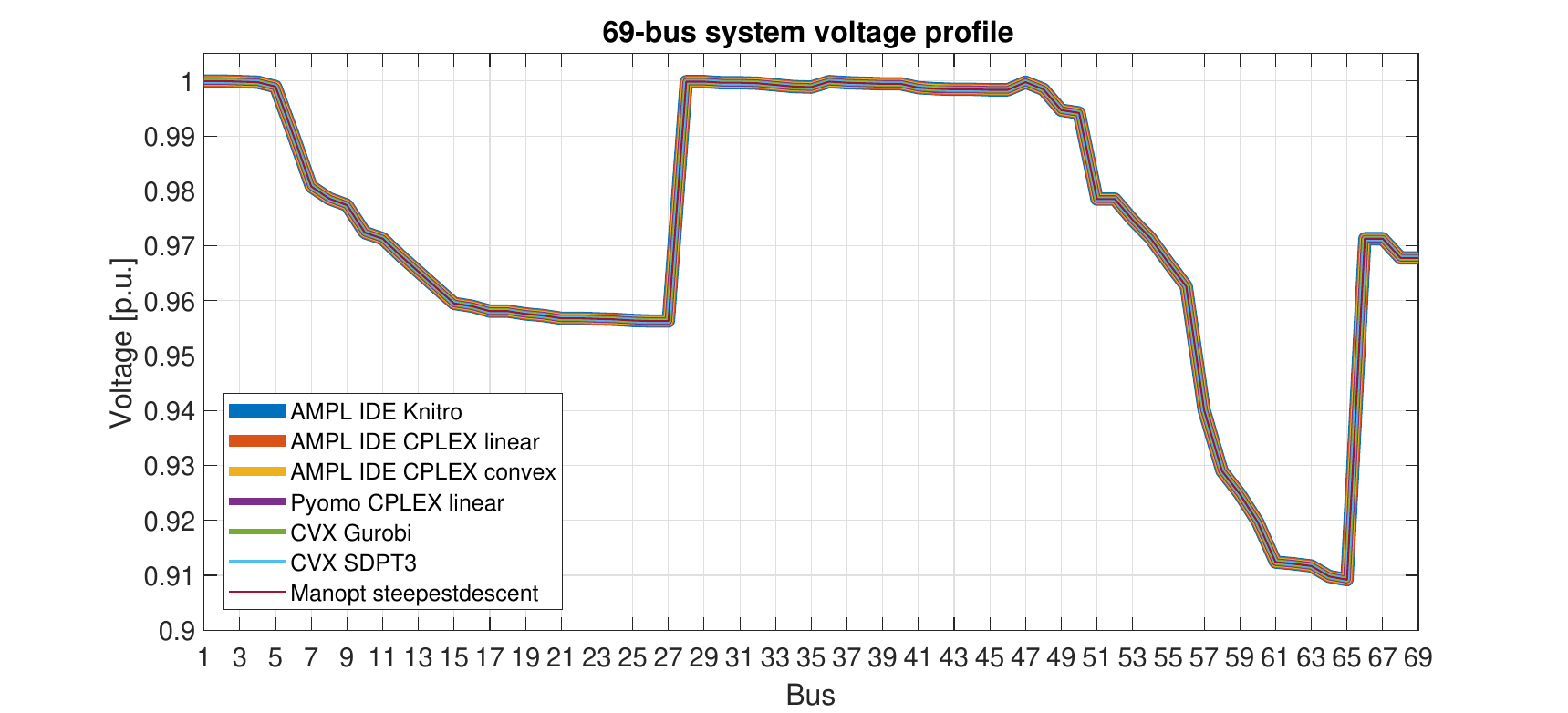}
\caption{Voltage profile for the 69-bus \gls{epds} power flow models.}
\label{FP:fig69busvolt}
\end{center}
\end{figure}

\begin{figure*}
\begin{center}
\includegraphics[trim=86pt 0pt 85pt 25pt,clip,width=.85\linewidth]{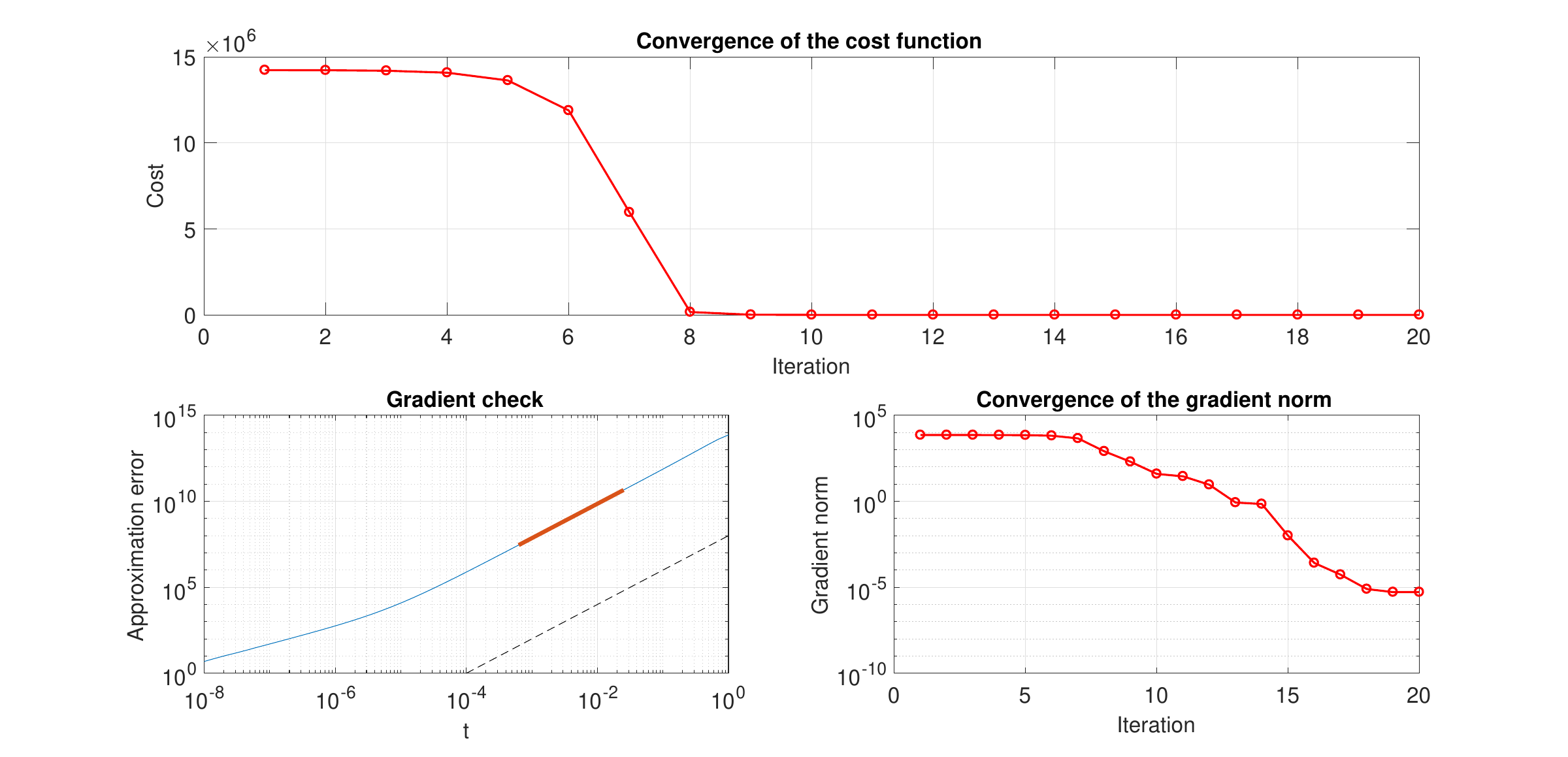}
\caption{Manopt metrics for the 69-bus \gls{epds} power flow \gls{bfs} model.}
\label{FP:manopt69}
\end{center}
\end{figure*}

Furthermore, Figure~\ref{FP:fig69busvolt} shows the voltage profiles. Since all the profiles overlap, it is again safe to assume that the models were well-formulated and well-conditioned. Also, Figure~\ref{FP:manopt69} presents the same Manopt metrics, showing the accuracy of the gradient, as well as the convergence of the cost function after 20 iterations and the gradient norm.

Thus, the gradient provided is accurate. Once again, the only models that exactly matched the value of the Knitro nonlinear model were Gurobi with the convex model and Manopt. Manopt once more achieved better results than CVX with SDPT3 and CPLEX in terms of accuracy. Therefore, Manopt has proven to be a powerful and reliable tool for solving the problem, not only matching other solvers, but also standing out for its robustness, flexibility, and accessibility in challenging scenarios. Its ability to deliver consistent results with very good execution time makes it a good choice for this power flow model in \gls{epds}. Additionally, Manopt has the potential to incorporate other planning strategies, such as reconfiguration, further expanding its applicability in power system optimization.

\subsection{\bf POWER FLOW IN ELECTRIC POWER TRANSMISSION SYSTEMS}

Power flow in \gls{epts} involves the analysis and determination of electrical quantities, such as voltage, current, and power, at different points in the network, ensuring that the generated energy is transmitted efficiently and safely. One of the main objectives of electric power systems is to meet the variable demand of consumers, providing energy with quality, i.e., with constant frequency and voltage. To achieve this, it is necessary to mathematically model the system's behavior, considering the characteristics of its components, such as generators, transformers, transmission lines, and loads. These elements are represented by nonlinear equations, which describe the relationships between electrical variables.

The analysis of power flow involves solving a system of nonlinear equations, which can be addressed using numerical methods such as the Newton-Raphson method, the decoupled method, the fast decoupled method, and others. These iterative methods allow for the determination of the voltages at the system's buses (magnitude and angle) that satisfy the active and reactive power balance equations \cite{romero}.

\begin{figure}
\begin{center}
\includegraphics[width=6cm]{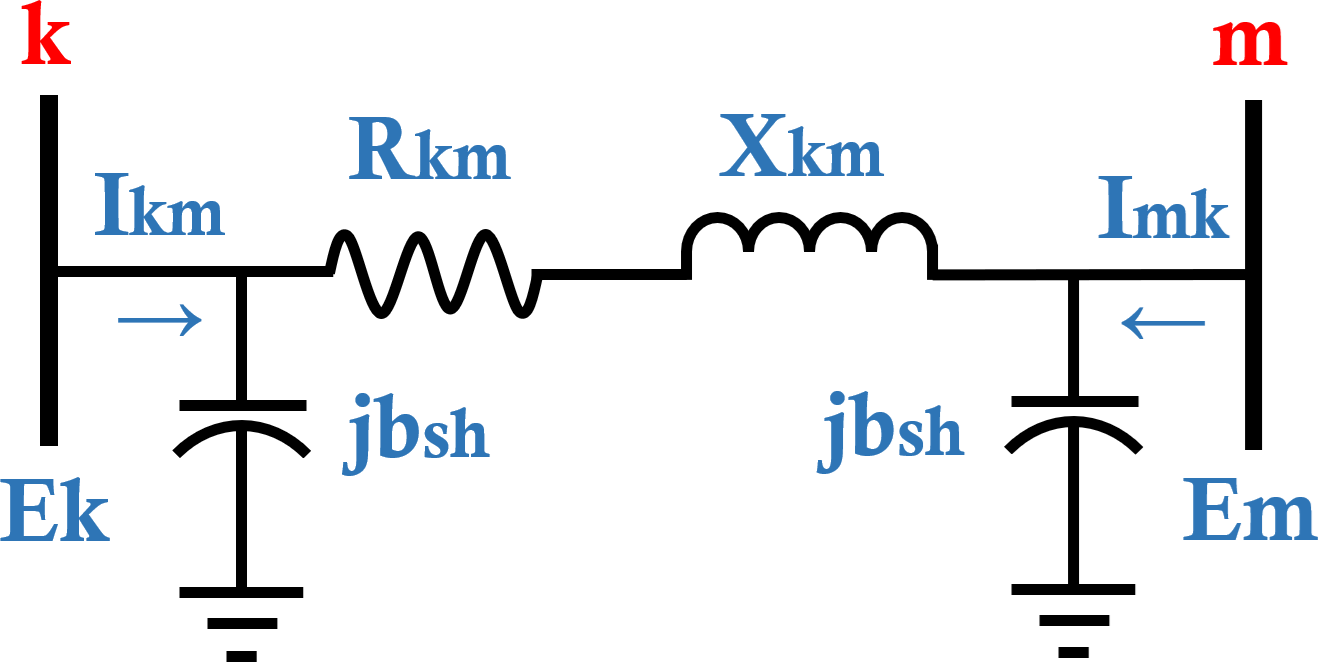}
\caption{Transmission line $\pi$-model.}
\label{pimodel}
\end{center}
\end{figure}

One of the most widely used mathematical models to represent a transmission line is the $\pi$-model, as suggested by Figure~\ref{pimodel}. In this case, the line's total capacitance is divided into two parts, placed at both the sending and receiving ends of the line, while the resistance and inductive reactance remain in the middle section of the line.

Thus, the active and reactive power balance equations are:
\begin{small}
\begin{align}
P_k &= V_k \sum_{m \in k} V_m (G_{km} \cos \theta_{km} + B_{km} \sin \theta_{km}), \label{Pk} \\
Q_k &= V_k \sum_{m \in k} V_m (G_{km} \sin \theta_{km} - B_{km} \cos \theta_{km}), \label{Qk}
\end{align}
\end{small}

\noindent where $P_k$ and $Q_k$ are the active and reactive power at bus $k$, $V_k$ and $V_m$ are the voltages at buses $k$ and $m$, $\theta_{km}$ is the angle difference between buses $k$ and $m$, and $G_{km}$ and $B_{km}$ are the conductance and susceptance of the line between buses $k$ and $m$.

The most commonly used numerical methods cited tend to linearize these nonlinear equations through a Taylor series expansion, generating a system of linear equations that can be solved iteratively. At each iteration, the methods update the state variables (voltages and angles) until the differences between the calculated values and the specified values (mismatches) are smaller than a predefined tolerance. It may require a good initial guess to ensure convergence, especially in systems with high load or more pronounced nonlinear characteristics. Additionally, the need to calculate and invert the Jacobian matrix at each iteration can be computationally expensive for very large systems \cite{romero2}.

In this context, \gls{mo} techniques can be applied to solve the power flow problem, since the power flow equations can be arranged in such a way that they are treated as equality constraints. Moreover, manifold-based methods can generalize the search for solutions, including the possibility of finding multiple or even all feasible solutions. This is because these methods are not limited to local convergence properties, as is often the case with Newton-Raphson. Furthermore, since p.u. values are generally used for calculations, the starting point can usually be an array of zeros, without the need to compute a precise starting point.

\subsubsection{\bf MANIFOLD SELECTION AND RIEMANNIAN GRADIENT PROJECTION}

It is desired to use the power flow equations to determine the variables of the \gls{epts} in Figure~\ref{FP:fig1}. Considering $S_{base} = 100$ MW as the base power, the bus and line data can be recovered from Tables~\ref{FP:31bus}-\ref{FP:31line}.

\begin{figure}
\begin{center}
\includegraphics[width=8cm]{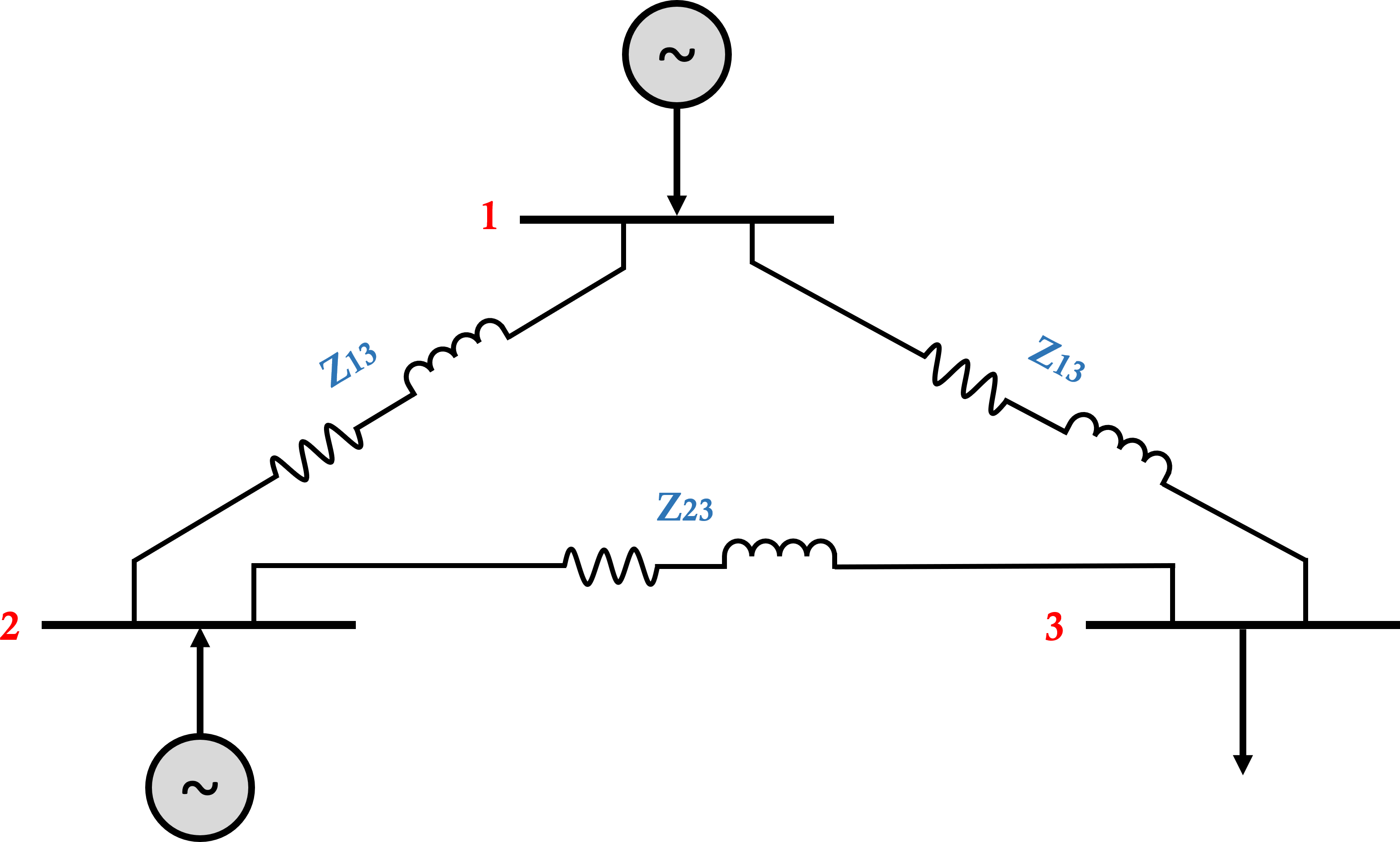}
\caption{1$^\circ$ 3-bus \gls{epts}.}
\label{FP:fig13bussys}
\end{center}
\end{figure}

\begin{table}
\centering
\begin{tabular}{ccccc}
\hline \hline
\textbf{Line} & \textbf{R (p.u)} & \textbf{X (p.u)} & \boldmath{$b_{shunt}/2$} \textbf{(p.u)} \\
\hline
1-2 & 0.01 & 0.03 & - \\
1-3 & 0.02 & 0.04 & - \\
2-3 & 0.0125 & 0.025 & - \\
\hline \hline
\end{tabular}
\caption{1$^\circ$ 3-bus \gls{epts} line data.}
\label{FP:31bus}
\end{table}

\begin{table*}
\centering
\resizebox{1.6\columnwidth}{!}{%
\begin{tabular}{ccccccccc}
\hline \hline
\textbf{Bus} & \textbf{Type} & \boldmath{$P_G$} \textbf{(MW)} & \boldmath{$Q_G$} \textbf{(MVAr)} & \boldmath{$P_D$} \textbf{(MW)} & \boldmath{$Q_D$} \textbf{(MVAr)} & \textbf{V (p.u)} & \boldmath{$\theta$} \textbf{(rad)} & \boldmath{$b_k$} \textbf{(p.u)} \\
\hline
1 & Slack & - & - & - & - & 1.05 & 0 & - \\
2 & PV & 200 & - & - & - & 1.04 & - & - \\
3 & PQ & - & - & 400 & 250 & - & - & - \\
\hline \hline
\end{tabular}}
\caption{1$^\circ$ 3-bus \gls{epts} bus data.}
\label{FP:31line}
\end{table*}

According to the Newton-Raphson method formulation, the Jacobian matrix for this system, aiming to solve subsystem 1 with respect to the unknown voltage magnitudes and angles, is:
\begin{small}
\begin{align}
\begin{bmatrix}
    \Delta P_2 \\
    \Delta P_3 \\
    \Delta Q_3
\end{bmatrix} =
\begin{bmatrix}
    H_{22} & H_{23} & N_{23} \\
    H_{32} & H_{33} & N_{33} \\
    M_{32} & M_{33} & L_{33}
\end{bmatrix}
\begin{bmatrix}
    \Delta \theta_2 \\
    \Delta \theta_3 \\
    \Delta V_3
\end{bmatrix}
.
\end{align}
\end{small}

Therefore, we have the equations for the differences between expected and calculated active and reactive power, which should approach zero:
\begin{small}
\begin{align}
    \Delta P_2 &= P_2^{exp} - P_2^{calc} = 2 - P_2^{calc}, \label{dP2} \\
    \Delta P_3 &= P_3^{exp} - P_3^{calc} = -4 - P_3^{calc}, \label{dP3} \\
    \Delta Q_3 &= Q_3^{exp} - Q_3^{calc} = -2.5 - Q_3^{calc} \label{dQ3},
\end{align}
\end{small}

\noindent where $P_2^{calc}$, $P_3^{calc}$ and $Q_3^{calc}$ are p.u. values obtained from Equations~\eqref{Pk}-\eqref{Qk}.

It is clear, therefore, that we have 6 variables: the subsystem 1 variables: $\theta_2$, $\theta_3$ and $V_3$, and later the subsystem 2 variables: $P_1^{calc}$, $Q_1^{calc}$ and $Q_2^{calc}$, which can be calculated using Equations~\eqref{Pk}-\eqref{Qk} after subsystem 1 has been solved.

In this way, we can use Equations~\eqref{Pk}-\eqref{Qk} in Equations~\eqref{dP2}-\eqref{dQ3} and set them equal to zero, and also set Equations~\eqref{Pk}-\eqref{Qk} for the subsystem 2 equal to zero to form a system of equations that refer to equality constraints. In this case, the cost function also has the form of Equation~\eqref{fcost}.

Since the Euclidean space can still be used for the Manifold, the optimization occurs in the space $\mathbb{R}^{2n}$, where $n$ corresponds to the number of buses in the \gls{epts}. The state of the power system is given by:
\begin{small}
\begin{align}
x = [\theta_2,~\theta_3,~V_3,~P_1^{calc},~Q_1^{calc},~Q_2^{calc}] \in \mathbb{R}^{2n},
\end{align}
\end{small}
and in a more generic notation, we have:
\begin{small}
\begin{align}
x = [\theta_{PV},~\theta_{PQ},~V_{PQ},~P_{V\theta}^{calc},~Q_{V\theta}^{calc},~Q_{PV}^{calc}] \in \mathbb{R}^{2n}.
\end{align}
\end{small}

We may hence define the following set:
\begin{small}
\begin{align}
\mathcal{M} = \{ x \in \mathbb{R}^{2n} \mid h(x) = 0 \}, \label{manifold3}
\end{align}
\end{small}

\noindent where $h(x)$ encodes the equations.

Thus, the projection of the gradient from the Euclidean space to the Riemannian space is trivial. By differentiating Equation~\eqref{fcost} considering that $h_i(x)$ are no longer complex-valued functions, we have:
\begin{small}
\begin{align}
\frac{\partial}{\partial x} |h_i(x)|^2 = 2 h_i(x) \nabla h_i(x)
\end{align}
\end{small}

Therefore, the gradient is:
\begin{small}
\begin{align}
\nabla \sum_{i=1}^{l+1} h_i(x)^2 =
\sum_{i=1}^{l+1} 2 h_i(x) \nabla h_i(x),
\end{align}
\end{small}

\noindent where $\nabla h_i(x)$ represents the $i$-th row of the Jacobian $J_h(x)$. In terms of the Jacobian:
\begin{small}
\begin{align}
\nabla \text{cost} = 2 J_h(x)^T h_i(x) \label{FP:grad}.
\end{align}
\end{small}

\subsubsection{\bf RESULTS WITH THE MANOPT TOOLBOX}

Using the 3-bus \gls{epts} shown in Figure~\ref{FP:fig13bussys} and the bus and line data that can be recovered from Tables~\ref{FP:31bus}-\ref{FP:31line}, the power flow can be solved using Manopt.

Therefore, we also want to compare the solution with the solution of commonly used numerical and approximate methods. Table~\ref{FP:tab3_1} presents the bus results for each method considering a tolerance of $10^{-3}$ for iterative methods.

\begin{table*}
\centering
\resizebox{1.75\columnwidth}{!}{%
\begin{tabular}{ccccc}
\hline \hline
    \textbf{Method} & \textbf{Bus voltage (p.u.)} & \textbf{Angle magnitudes (rad)} & \textbf{Bus power (p.u.)} & \textbf{Iterations} \\
    \hline \hline
    & $V_1 = 1.0500$ & $|\theta_1| = 0.0000$ & $P_1 = 2.1840$ $Q_1 = 1.4085$ & \\
    \textbf{Newton-Raphson} & $V_2 = 1.0400$ & $|\theta_2| = 0.0087$ & $P_2 = 2$ $Q_2 = 1.4616$ & 3 \\
    & $V_3 = 0.9717$ & $|\theta_3| = 0.0471$ & $P_3 = -4$ $Q_3 = -2.5$ & \\ \hline
    & $V_1 = 1.0500$ & $|\theta_1| = 0.0000$ & $P_1 = 2.1812$ $Q_1 = 1.4078$ & \\
    \textbf{Decoupled} & $V_2 = 1.0400$ & $|\theta_2| = 0.0087$ & $P_2 = 2$ $Q_2 = 1.4613$ & 5 P$\theta$ and 5 QV \\
    & $V_3 = 0.9718$ & $|\theta_3| = 0.0469$ & $P_3 = -4$ $Q_3 = -2.5$ & \\ \hline
    & $V_1 = 1.0500$ & $|\theta_1| = 0.0000$ & $P_1 = 2.1835$ $Q_1 = 1.4091$ & \\
    \textbf{Fast Decoupled} & $V_2 = 1.0400$ & $|\theta_2| = 0.0087$ & $P_2 = 2$ $Q_2 = 1.4619$ & 4 P$\theta$ and 4 QV \\
    & $V_3 = 0.9717$ & $|\theta_3| = 0.0471$ & $P_3 = -4$ $Q_3 = -2.5$ & \\ \hline
    & & $|\theta_1| = 0.0000$ & $P_1 = 2.0000$ & \\
    \textbf{DC} & - & $|\theta_2| = 0.0095$ & $P_2 = 2.0000$ & - \\
    & & $|\theta_3| = 0.0674$ & $P_3 = -4$ & \\ \hline
    & & $|\theta_1| = 0.0000$ & $P_1 = 2.0536$ & \\
    \textbf{DC with losses} & - & $|\theta_2| = 0.0107$ & $P_2 = 1.9727$ & - \\
    & & $|\theta_3| = 0.0689$ & $P_3 = -4.0495$ & \\ \hline
    \textbf{Manopt} & $V_1 = 1.0500$ & $|\theta_1| = 0.0000$ & $P_1 = 2.1842$ $Q_1 = 1.4085$ & \\
    \textbf{(trustregions solver)} & $V_2 = 1.0400$ & $|\theta_2| = 0.0087$ & $P_2 = 2$ $Q_2 = 1.4618$ & 14 \\
    & $V_3 = 0.9717$ & $|\theta_3| = 0.0471$ & $P_3 = -4$ $Q_3 = -2.5$ & \\ \hline \hline
\end{tabular}}
\caption{1$^\circ$ 3-bus \gls{epts} results for each method.}
\label{FP:tab3_1}
\end{table*}

\begin{figure*}
\begin{center}
\includegraphics[trim=86pt 0pt 85pt 25pt,clip,width=.9\linewidth]{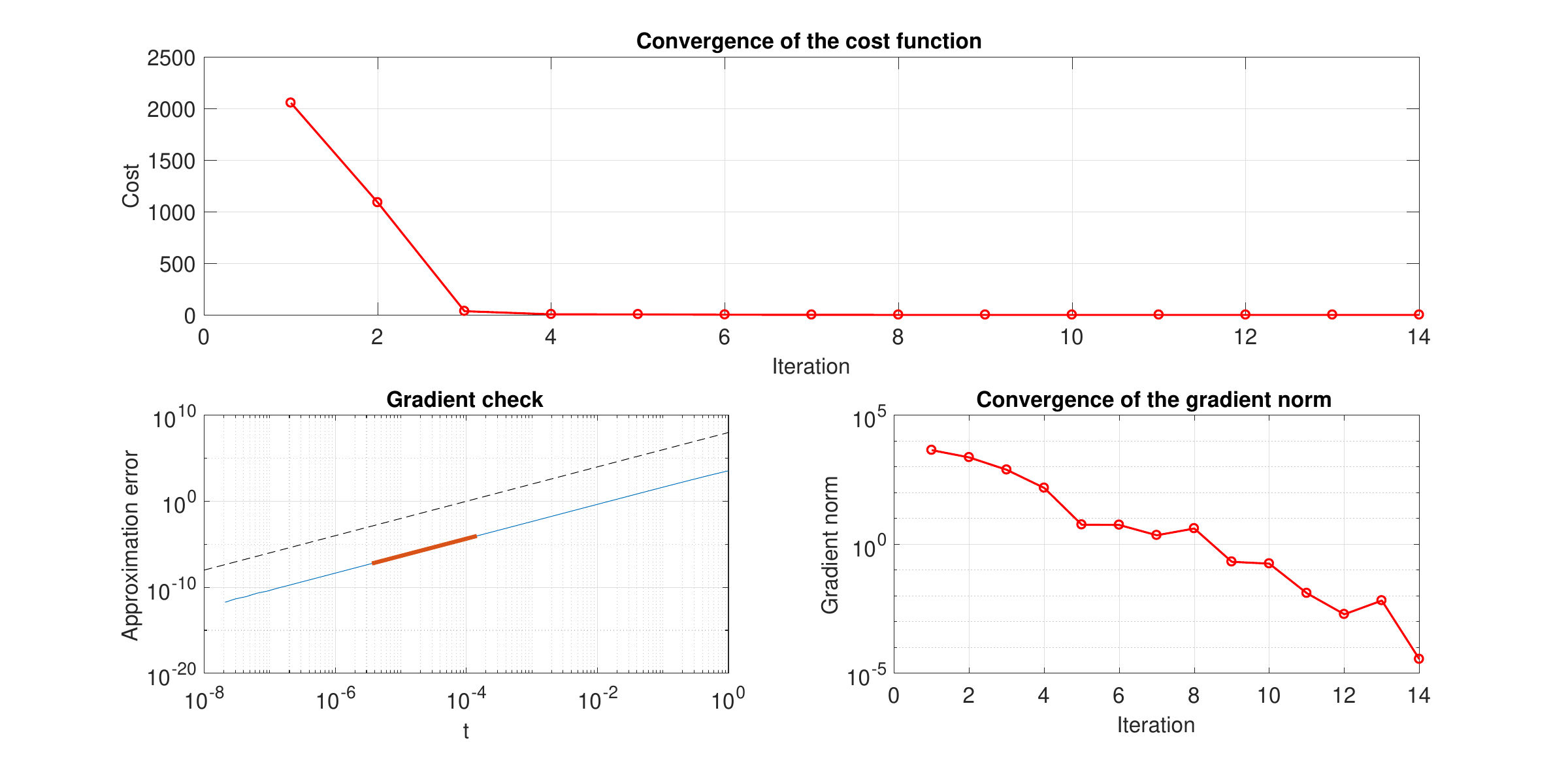}
\caption{Manopt metrics for the 1$^\circ$ 3-bus \gls{epts} power flow approach.}
\label{FP:manopt31}
\end{center}
\end{figure*}

Furthermore, Figure~\ref{FP:manopt31} presents the Manopt metrics for the 3-bus \gls{epts} power flow approach, showing the accuracy of the gradient and the Hessian, as well as the convergence of the cost function in 14 iterations and the gradient norm.

It can be concluded that Manopt’s trustregions solver delivered excellent performance, reaching a cost function value of $3.858995\cdot 10^{-10}$ in 0.199434 seconds, a convergence behavior akin to the Newton-Raphson method, as desired. Although computational time is not being evaluated here, since the DC method would certainly outperform all of them because it is not iterative, it is important to note that Manopt continues to have a very low computational time, which is great if the idea is to use this approach and expand it to systems planning.

\begin{figure}
\begin{center}
\includegraphics[width=8cm]{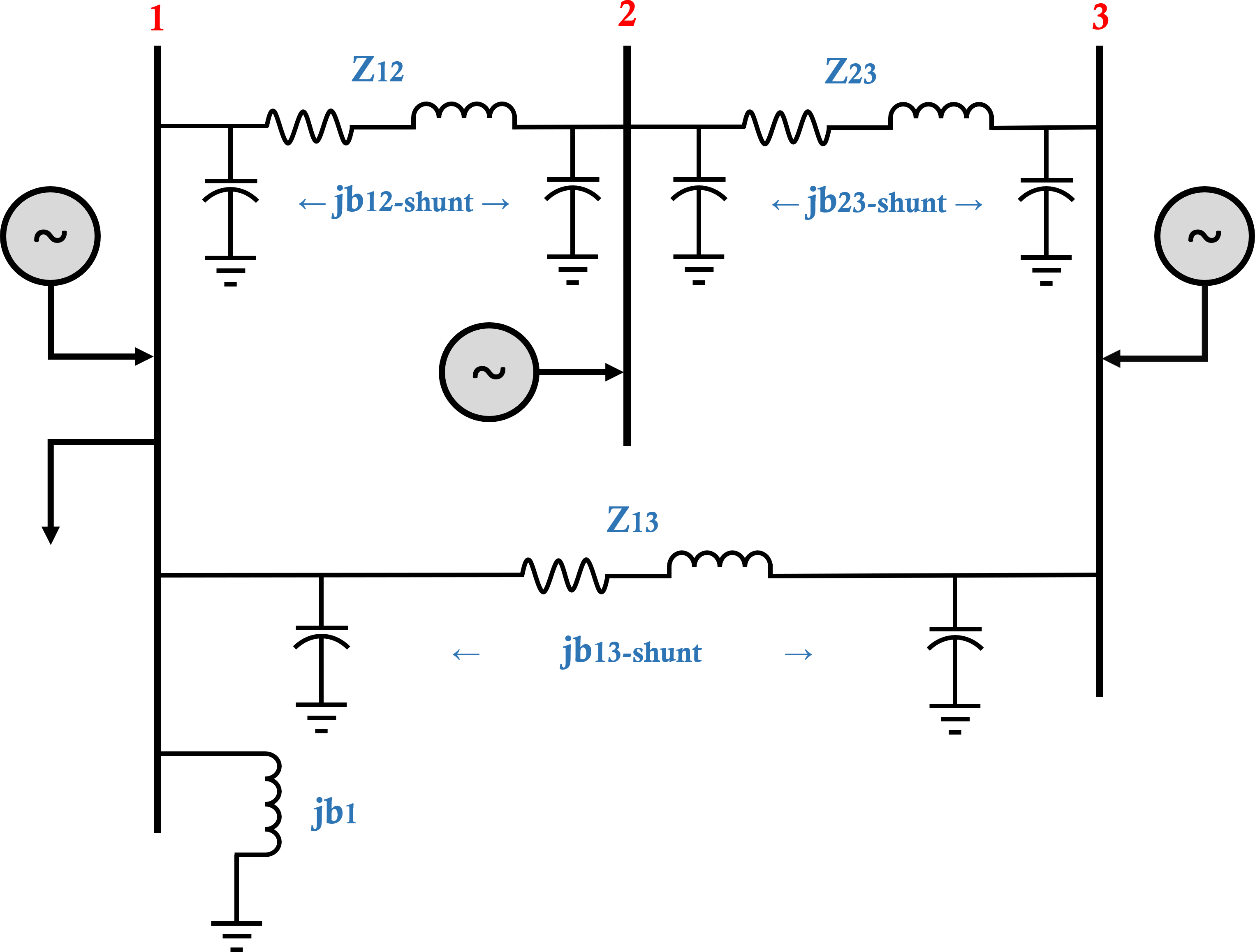}
\caption{2$^\circ$ 3-bus \gls{epts}.}
\label{FP:fig23bussys}
\end{center}
\end{figure}

Subsequently, in order to verify another test case, using the 3-bus \gls{epts} shown in Figure~\ref{FP:fig23bussys} and the bus and line data that can be recovered from Tables~\ref{FP:32bus}-\ref{FP:32line}, the power flow can be solved using Manopt.

\begin{table}
\centering
\begin{tabular}{ccccc}
\hline \hline
\textbf{Line} & \textbf{R (p.u)} & \textbf{X (p.u)} & \boldmath{$b_{shunt}/2$} \textbf{(p.u)} \\
\hline
1-2 & 0.03 & 0.3 & j0.02 \\
1-3 & 0.08 & 1.1 & j0.03 \\
2-3 & 0.05 & 0.8 & j0.01 \\
\hline \hline
\end{tabular}
\caption{2$^\circ$ 3-bus \gls{epts} line data.}
\label{FP:32bus}
\end{table}

\begin{table*}
\centering
\resizebox{1.4\columnwidth}{!}{%
\begin{tabular}{ccccccccc}
\hline \hline
\textbf{Bus} & \textbf{Type} & \boldmath{$P_G$} \textbf{(MW)} & \boldmath{$Q_G$} \textbf{(MVAr)} & \boldmath{$P_D$} \textbf{(MW)} & \boldmath{$Q_D$} \textbf{(MVAr)} & \textbf{V (p.u)} & \boldmath{$\theta$} \textbf{(rad)} & \boldmath{$b_k$} \textbf{(p.u)} \\
\hline
1 & PQ & - & 5 & 30 & - & - & - & -j0.05 \\
2 & Slack & - & - & - & - & 1.00 & 0 & - \\
3 & PV & 20 & - & - & - & 1.00 & - & - \\
\hline \hline
\end{tabular}}
\caption{2$^\circ$ 3-bus \gls{epts} bus data.}
\label{FP:32line}
\end{table*}

Table~\ref{FP:tab3_2} presents the bus results for each method considering a tolerance of $10^{-3}$ for iterative methods.

\begin{table*}
\centering
\resizebox{1.74\columnwidth}{!}{%
\begin{tabular}{ccccc}
\hline \hline
\textbf{Method} & \textbf{Bus voltage (p.u.)} & \textbf{Angle magnitudes (rad)} & \textbf{Bus power (p.u.)} & \textbf{Iterations} \\
\hline \hline
& $V_1 = 1.0025$ & $|\theta_1| = 0.0569$ & $P_1 = -0.3$ $Q_1 = 0.05$ & \\
\textbf{Newton-Raphson} & $V_2 = 1.0000$ & $|\theta_2| = 0.0000$ & $P_2 = 0.1025$ $Q_2 = -0.0433$ & 3 \\
& $V_3 = 1.0000$ & $|\theta_3| = 0.0689$ & $P_3 = 0.2$ $Q_3 = -0.0457$ & \\ \hline
& $V_1 = 1.0025$ & $|\theta_1| = 0.0564$ & $P_1 = -0.3$ $Q_1 = 0.05$ & \\
\textbf{Decoupled} & $V_2 = 1.0000$ & $|\theta_2| = 0.0000$ & $P_2 = 0.1003$ $Q_2 = -0.0434$ & 2 P$\theta$ and 2 QV \\
& $V_3 = 1.0000$ & $|\theta_3| = 0.0693$ & $P_3 = 0.2$ $Q_3 = -0.0458$ & \\ \hline
& $V_1 = 1.0027$ & $|\theta_1| = 0.0559$ & $P_1 = -0.3$ $Q_1 = 0.05$ & \\
\textbf{Fast Decoupled} & $V_2 = 1.0000$ & $|\theta_2| = 0.0000$ & $P_2 = 0.0988$ $Q_2 = -0.0437$ & 2 P$\theta$ and 2 QV \\
& $V_3 = 1.0000$ & $|\theta_3| = 0.0691$ & $P_3 = 0.2$ $Q_3 = -0.0459$ & \\ \hline
& & $|\theta_1| = 0.0559$ & $P_1 = -0.3000$ & \\
\textbf{DC} & - & $|\theta_2| = 0.0000$ & $P_2 = 0.1000$ & - \\
& & $|\theta_3| = 0.0691$ & $P_3 = 0.2000$ & \\ \hline
& & $|\theta_1| = 0.0563$ & $P_1 = -0.3010$ & \\
\textbf{DC with losses} & - & $|\theta_2| = 0.0000$ & $P_2 = 0.1010$ & - \\
& & $|\theta_3| = 0.0686$ & $P_3 = 0.1993$ & \\ \hline
& $V_1 = 1.0025$ & $|\theta_1| = 0.0569$ & $P_1 = -0.3$ $Q_1 = 0.05$ & \\
\textbf{Manopt} & $V_2 = 1.0000$ & $|\theta_2| = 0.0000$ & $P_2 = 0.1024$ $Q_2 = -0.0433$ & 8 \\
\textbf{(trustregions solver)} & $V_3 = 1.0000$ & $|\theta_3| = 0.0689$ & $P_3 = 0.2$ $Q_3 = -0.0458$ & \\ \hline \hline
\end{tabular}}
\caption{2$^\circ$ 3-bus \gls{epts} results for each method.}
\label{FP:tab3_2}
\end{table*}

\begin{figure*}
\begin{center}
\includegraphics[trim=86pt 0pt 85pt 25pt,clip,width=.9\linewidth]{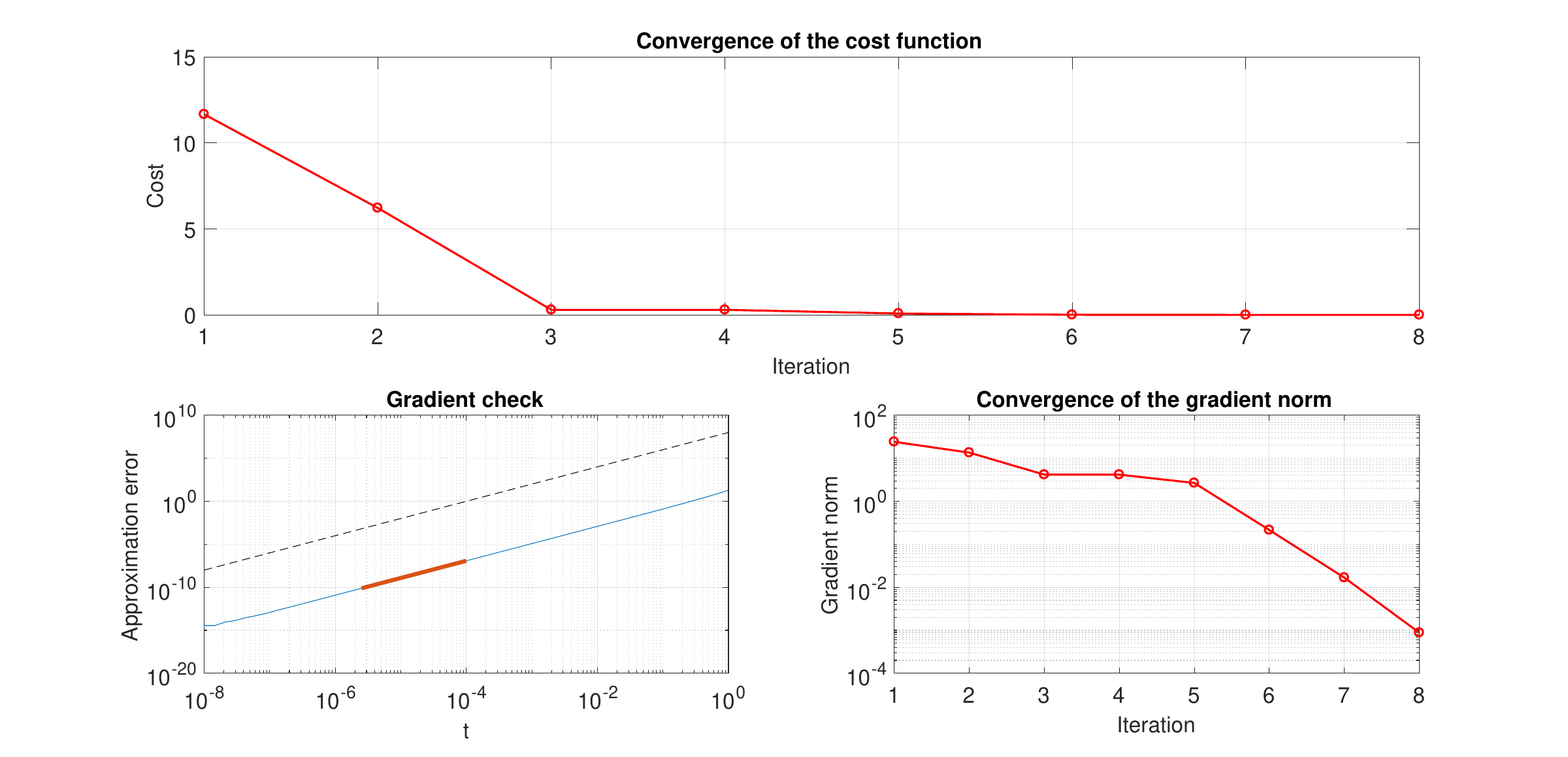}
\caption{Manopt metrics for the 2$^\circ$ 3-bus \gls{epts} power flow approach.}
\label{FP:manopt32}
\end{center}
\end{figure*}

Furthermore, Figure~\ref{FP:manopt32} presents the Manopt metrics for the 3-bus \gls{epts} power flow approach, showing the accuracy of the gradient and the Hessian, as well as the convergence of the cost function in 8 iterations and the gradient norm.

Once again, Manopt’s trustregions solver delivered excellent performance, reaching a cost function value of $7.609651\cdot 10^{-10}$ in 0.191431 seconds, a convergence behavior akin to the Newton-Raphson method, as desired.

\begin{figure}
\begin{center}
\includegraphics[width=8.5cm]{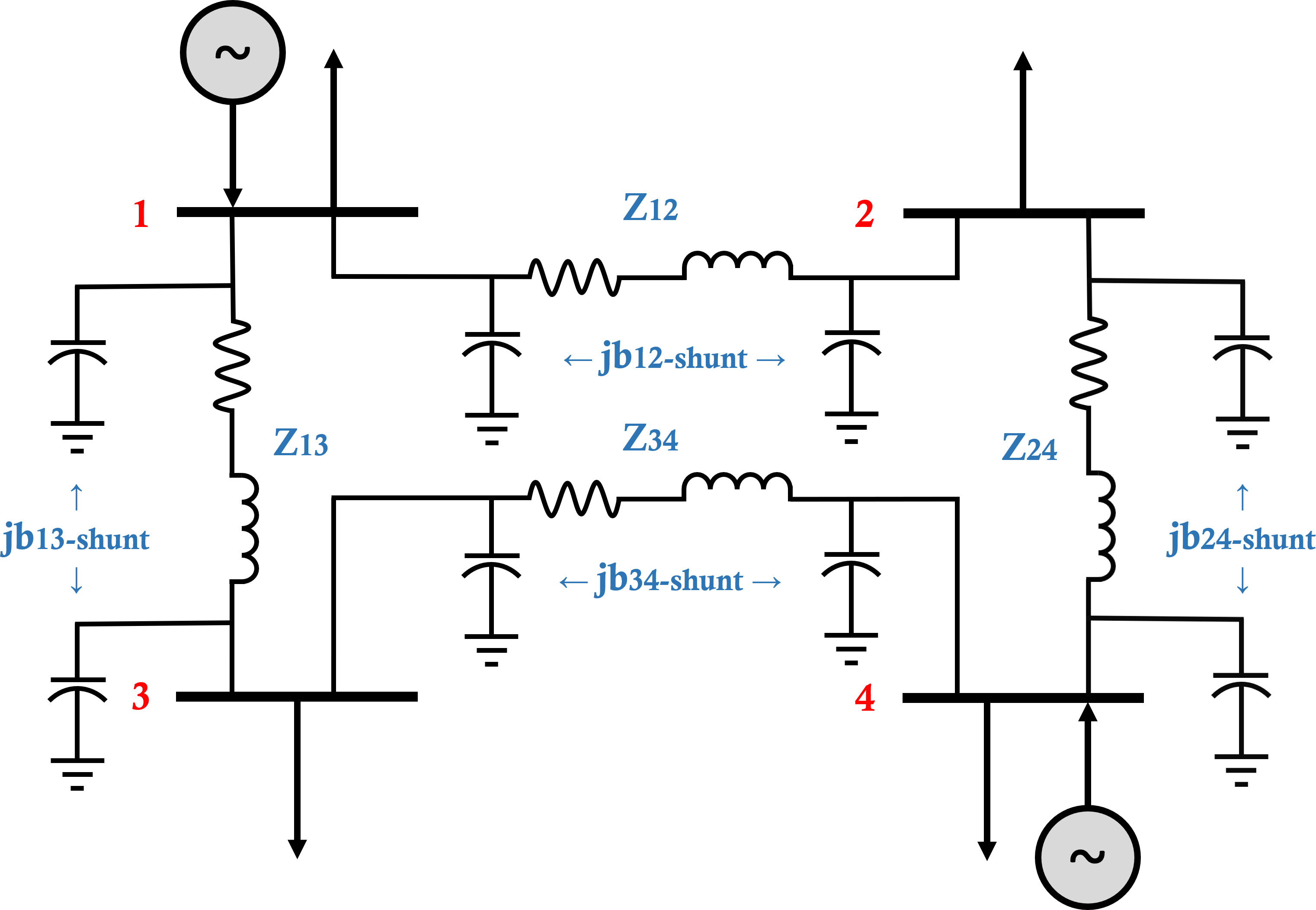}
\caption{4-bus \gls{epts}.}
\label{FP:fig4bussys}
\end{center}
\end{figure}

Finally, in order to verify one last test case, using the 4-bus \gls{epts} shown in Figure~\ref{FP:fig4bussys} and the bus and line data that can be recovered from Tables~\ref{FP:4bus}-\ref{FP:4line}, the power flow can be solved using Manopt.

\begin{table}
\centering
\begin{tabular}{ccccc}
\hline \hline
\textbf{Line} & \textbf{R (p.u)} & \textbf{X (p.u)} & \boldmath{$b_{shunt}$} \textbf{(p.u)} \\
\hline
1-2 & 0.01008 & 0.05040 & j0.05125 \\
1-3 & 0.00744 & 0.03720 & j0.03875 \\
2-4 & 0.00744 & 0.03720 & j0.03875 \\
3-4 & 0.01272 & 0.06360 & j0.06375 \\
\hline \hline
\end{tabular}
\caption{4-bus \gls{epts} line data.}
\label{FP:4bus}
\end{table}

\begin{table*}
\centering
\resizebox{1.55\columnwidth}{!}{%
\begin{tabular}{ccccccccc}
\hline \hline
\textbf{Bus} & \textbf{Type} & \boldmath{$P_G$} \textbf{(MW)} & \boldmath{$Q_G$} \textbf{(MVAr)} & \boldmath{$P_D$} \textbf{(MW)} & \boldmath{$Q_D$} \textbf{(MVAr)} & \textbf{V (p.u)} & \boldmath{$\theta$} \textbf{(rad)} & \boldmath{$b_k$} \textbf{(p.u)} \\
\hline
1 & Slack & - & - & 50 & 30.99 & 1.00 & 0 & - \\
2 & PQ & 0 & 0 & 170 & 105.35 & 1.00 & 0 & - \\
3 & PQ & 0 & 0 & 200 & 123.94 & 1.00 & 0 & - \\
4 & PV & 318 & - & 80 & 49.58 & 1.02 & 0 & - \\
\hline \hline
\end{tabular}}
\caption{4-bus \gls{epts} bus data.}
\label{FP:4line}
\end{table*}

Table~\ref{FP:tab4_1} presents the bus results for each method considering a tolerance of $10^{-3}$ for iterative methods.

\begin{table*}
\centering
\resizebox{1.8\columnwidth}{!}{%
\begin{tabular}{ccccc}
\hline \hline
\textbf{Method} & \textbf{Bus voltage (p.u.)} & \textbf{Angle magnitudes (rad)} & \textbf{Bus power (p.u.)} & \textbf{Iterations} \\
\hline \hline
& $V_1 = 1.0000$ & $|\theta_1| = 0.0000$ & $P_1 = 1.3679$ $Q_1 = 0.8350$ & \\
\textbf{Newton-Raphson} & $V_2 = 0.9824$ & $|\theta_2| = 0.0170$ & $P_2 = -1.7$ $Q_2 = -1.0535$ & 3 \\
& $V_3 = 0.9690$ & $|\theta_3| = 0.0327$ & $P_3 = -2.0$ $Q_3 = -1.2394$ & \\
& $V_4 = 1.0200$ & $|\theta_4| = 0.0266$ & $P_4 = 2.38$ $Q_4 = 1.3184$ & \\ \hline
& $V_1 = 1.0000$ & $|\theta_1| = 0.0000$ & $P_1 = 1.3665$ $Q_1 = 0.8350$ & \\
\textbf{Decoupled} & $V_2 = 0.9824$ & $|\theta_2| = 0.0170$ & $P_2 = -1.7$ $Q_2 = -1.0535$ & 3 P$\theta$ and 3 QV \\
& $V_3 = 0.9690$ & $|\theta_3| = 0.0327$ & $P_3 = -2.0$ $Q_3 = -1.2394$ & \\
& $V_4 = 1.0200$ & $|\theta_4| = 0.0266$ & $P_4 = 2.38$ $Q_4 = 1.3183$ & \\ \hline
& $V_1 = 1.0000$ & $|\theta_1| = 0.0000$ & $P_1 = 1.3740$ $Q_1 = 0.8353$ & \\
\textbf{Fast Decoupled} & $V_2 = 0.9824$ & $|\theta_2| = 0.0171$ & $P_2 = -1.7$ $Q_2 = -1.0535$ & 3 P$\theta$ and 3 QV \\
& $V_3 = 0.9690$ & $|\theta_3| = 0.0328$ & $P_3 = -2.0$ $Q_3 = -1.2394$ & \\
& $V_4 = 1.0200$ & $|\theta_4| = 0.0266$ & $P_4 = 2.38$ $Q_4 = 1.3189$ & \\ \hline
& & $|\theta_1| = 0.0000$ & $P_1 = 1.3200$ & \\
\textbf{DC} & - & $|\theta_2| = 0.0185$ & $P_2 = -1.7000$ & - \\
& & $|\theta_3| = 0.0355$ & $P_3 = -2.0000$ & \\
& & $|\theta_4| = 0.0311$ & $P_4 = 2.3800$ & \\ \hline
& & $|\theta_1| = 0.0000$ & $P_1 = 1.3461$ & \\
\textbf{DC with losses} & - & $|\theta_2| = 0.0192$ & $P_2 = -1.7070$ & - \\
& & $|\theta_3| = 0.0361$ & $P_3 = -2.0100$ & \\
& & $|\theta_4| = 0.0301$ & $P_4 = 2.3669$ & \\ \hline
& $V_1 = 1.0000$ & $|\theta_1| = 0.0000$ & $P_1 = 1.3681$ $Q_1 = 0.8351$ & \\
\textbf{Manopt} & $V_2 = 0.9824$ & $|\theta_2| = 0.0170$ & $P_2 = -1.7$ $Q_2 = -1.0535$ & 15 \\
\textbf{(trustregions solver)} & $V_3 = 0.9690$ & $|\theta_3| = 0.0327$ & $P_3 = -2.0$ $Q_3 = -1.2394$ & \\
& $V_4 = 1.0200$ & $|\theta_4| = 0.0266$ & $P_4 = 2.38$ $Q_4 = 1.3185$ & \\ \hline \hline
\end{tabular}}
\caption{4-bus \gls{epts} results for each method.}
\label{FP:tab4_1}
\end{table*}

\begin{figure*}
\begin{center}
\includegraphics[trim=86pt 0pt 85pt 25pt,clip,width=.85\linewidth]{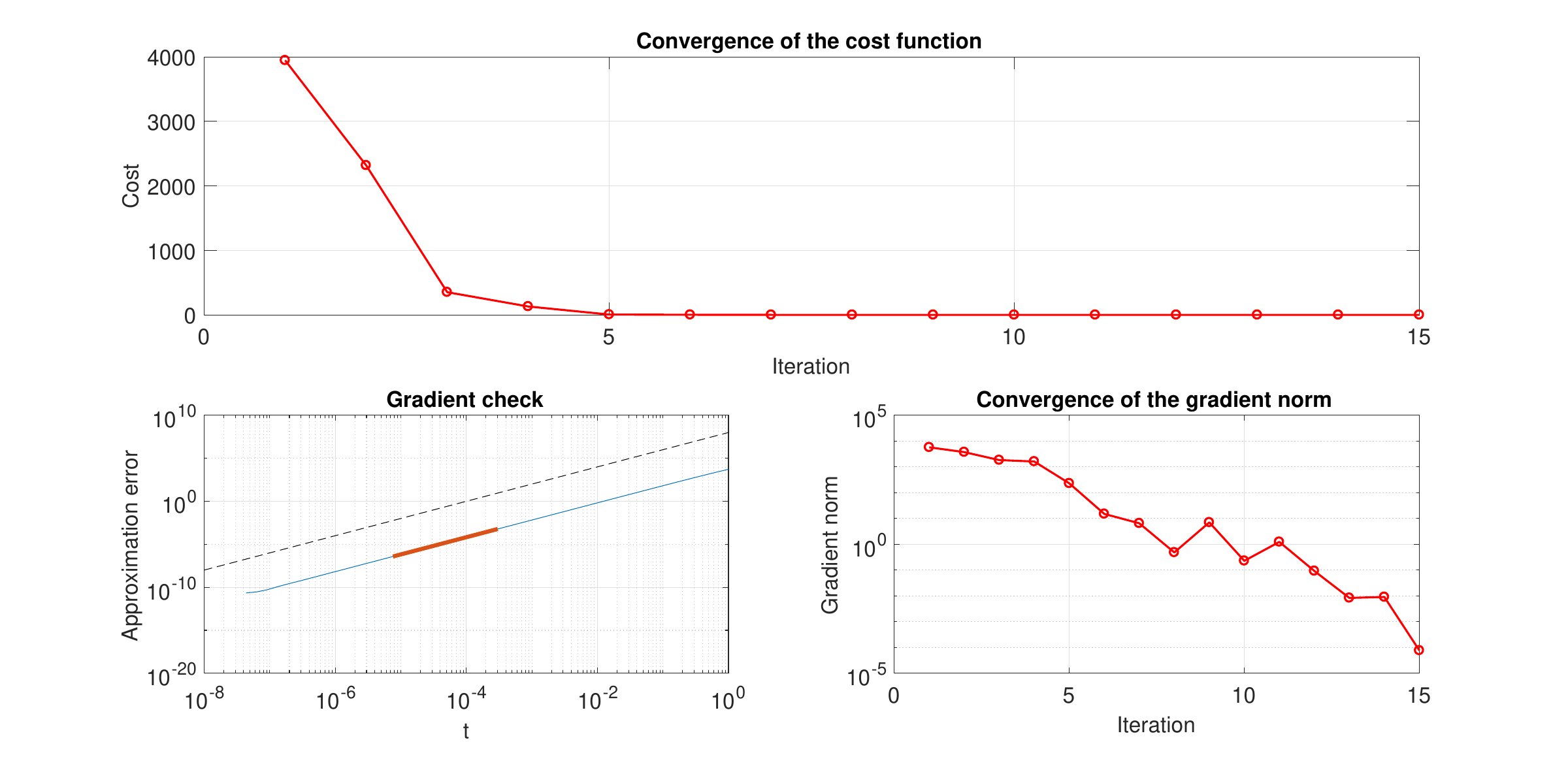}
\caption{Manopt metrics for the 4-bus \gls{epts} power flow approach.}
\label{FP:manopt4}
\end{center}
\end{figure*}

Furthermore, Figure~\ref{FP:manopt4} presents the Manopt metrics for the 4-bus \gls{epts} power flow approach, showing the accuracy of the gradient and the Hessian, as well as the convergence of the cost function in 15 iterations and the gradient norm.

Once again, Manopt’s trustregions solver delivered excellent performance, reaching a cost function value of $6.920804\cdot 10^{-10}$ in 0.202117 seconds, a convergence behavior akin to the Newton-Raphson method, as desired. These results reinforce Manopt as a powerful and reliable tool for optimization in \gls{epts}, combining numerical robustness with computational efficiency. Its flexibility in dealing with problem scalability positions it as an attractive alternative to traditional methods.

\section{INTEGRATION CHALLENGES OF MANIFOLD OPTIMIZATION TECHNIQUES IN POWER SYSTEMS}
\label{sec:BarriersMO}

While \gls{mo} offers significant advantages for power systems, its integration into existing infrastructure does present some challenges. This section addresses these challenges, drawing insights from the successful application of \gls{mo} in power system optimization presented in this paper.

\vspace{2mm}\noindent\textbf{{Complexity of Implementation}}. The theoretical foundations of \gls{mo} involve differential geometry and manifold theory, which can appear complex. At first, the classical mathematical model of power flow in \gls{epds} may seem very complex. Including all the restrictions in the cost function, deriving it on the chosen Manifold and computing an initial point for it may be challenging. Furthermore, deriving a model from the power equations in \gls{epts}, given the wide variety of traditional methods, can be challenging.\\
\textit{Mitigation:} The practical implementation of MO, as demonstrated in this paper for power flow problems, has proven to be less intricate than initially perceived. The use of tools like Manopt, reducing the need for deep specialized knowledge in differential geometry, with the choice of using the \gls{bfs} approach for \gls{epds}, simplifies the process. The presented power flow approaches for the \gls{epds} and \gls{epts} were all successfully solved using Manopt, highlighting a relatively straightforward implementation process once the problem is formulated.

\vspace{2mm}
\noindent\textbf{{Algorithmic Design and Tuning:}} The classical mathematical power flow model would undoubtedly pose challenges when incorporating all constraints into the cost function, deriving it from the chosen manifold, and computing an initial point. Furthermore, the very nature of the \gls{bfs} approach can be inherently challenging, considering that the problem variables are usually complex and must be treated in such a way as to separate the real and imaginary components, and also the possibility of this approach including expansion planning, which is generally not done in the literature.\\
\textit{Mitigation:} The case studies demonstrate that these algorithms converge efficiently, indicating a \textit{plug-and-play} capability for many power system scenarios once the problem is appropriately structured as an optimization on a manifold. This contradicts the notion that significant expertise in algorithm selection and tuning is always required.

\vspace{2mm}\noindent\textbf{{Computational Overhead}}
It stands to reason that for large-scale power systems, the Manopt implemented algorithm would exhibit scaled execution time. In any case, power flow problems tend to be inherently complex regardless of the approaches chosen, so papers focus heavily on computational time and accuracy.\\
\textit{Mitigation:} For the 14-bus, 33-bus, and 69-bus \gls{epds} models, Manopt achieved comparable or better execution times while maintaining very high accuracy. This suggests that \gls{mo} methods do not inherently lead to increased processing times or hinder real-time processing capabilities in power systems. The tendency is that the computational time is kept low for \gls{epds} and \gls{epts}, and is of the same order as other commercial solvers, maintaining competitiveness.

\vspace{2mm}\noindent\textbf{{Integration with Existing Systems}}
A legitimate concern would certainly be the use of standard IEEE test systems, given their complexity and large scale.\\
\textit{Mitigation:} The successful application of \gls{mo} to power flow problems in \gls{epds} and \gls{epts} within this paper demonstrates its potential for seamless integration with existing power system analysis frameworks. By formulating power flow equations as equality constraints on a manifold with the shown approaches, \gls{mo} tends to be able to solve any power system with almost perfect accuracy and with very good computational time.

\vspace{2mm}\noindent\textbf{{Scalability Concerns}}
The complexity of the optimization problem grows with the number of users in a network. \gls{mo} methods must be scalable to handle large-scale deployments without compromising performance. A common concern would certainly be the scalability of the algorithm when applied to large-scale, complex systems such as IEEE standard test systems.\\
\textit{Mitigation:} The scalability of \gls{mo} methods has been positively highlighted in this paper. The results for 14-bus, 33-bus, and 69-bus \gls{epds} models demonstrate that \gls{mo} performs consistently well across different system sizes, with computational times remaining highly competitive. This indicates that \gls{mo} can effectively handle the increasing complexity and dimensionality of modern power systems, positioning it as a robust solution for large-scale deployments.

\vspace{2mm}\noindent\textbf{{Planning for Expansion}}
The \gls{epds} and \gls{epts} planning problem is a classical challenge characterized by a combinatorial explosion, which has rightly attracted significant attention from researchers employing diverse techniques to tackle such complexity. The primary difficulties in solving this problem stem from the combinatorial nature of the planning process, which typically leads to an explosive number of alternatives, even for medium-sized systems \cite{romero}.\\
\textit{Mitigation:} Manopt has the potential to incorporate other planning strategies, such as optimal reconfiguration, further expanding its applicability in power system optimization. One approach could involve a Master-Slave configuration, where the Master proposes a value for the problem’s core variable (often binary or integer), and the Slave returns the corresponding power flow calculation. Considering Algorithm~\ref{alg:bfs_power_flow}, it would be enough to consider a vector of binary variables $y$, which represents the activation or inactivation of the lines, given by the Master and use $I_L(i) \gets I_n^{aux}(N_j(i)) \cdot y(i)$ in Algorithm~\ref{alg:bfs_power_flow}.\\
\textit{Data ordering:} It is said that to work correctly with the \gls{bfs} approach, the line data must be ordered, which is typically not the case with reconfiguration or other strategies. A simple algorithm with negligible computational overhead, such as Algorithm~\ref{alg:branch_reorganization}, can be used to organize the data in the data transition from Master to Slave.

\begin{algorithm}
\small
\caption{Branch reorganization}
\label{alg:branch_reorganization}
\begin{algorithmic}[1]
\State \textbf{Input:} 
\State \quad $L$ - List of all branches (sender and receiver)
\State \quad $sub$ - Slack bus/substation index

\State \textbf{Output:} 
\State \quad $L$ - Reorganized branch list

\State $sub_{aux} \gets \{sub\}$ \Comment{Initialize with substation}

\For{$k \gets 1$ \textbf{to} $\text{length}(L)$}
    \For{$i \gets 1$ \textbf{to} $\text{size}(L,1)$}
        \If{$L(i,2) \in sub_{aux}$ \textbf{and} $L(i,1) \notin sub_{aux}$}
            \State $L(i,1:2) \gets [L(i,2), L(i,1)]$ \Comment{Swap sender and receiver}
        \EndIf
        \If{$L(i,1) \in sub_{aux}$ \textbf{and} $L(i,2) \notin sub_{aux}$}
            \State $sub_{aux} \gets sub_{aux} \cup \{L(i,2)\}$ \Comment{Add new accessible bus}
        \EndIf
    \EndFor
\EndFor
\end{algorithmic}
\end{algorithm}

\section{CONCLUSION} \label{sec:VII}

This paper demonstrates that \gls{mo} offers a robust and versatile framework for tackling complex optimization challenges, particularly the non-convex power flow problems inherent in \gls{epds} and \gls{epts}. \gls{mo} has shown significant advantages over conventional methods by leveraging the inherent geometric properties of manifolds. Specifically, it reformulates constrained problems as unconstrained optimization over a smooth manifold structure. This approach, as applied in this work, facilitates more efficient and effective solutions by naturally handling non-convex constraints and exploiting these geometric properties, proving particularly well-suited for the intricacies of advanced power systems.

The application of \gls{mo} in \gls{epds} reveals that when applied to power flow optimization, demonstrates highly competitive performance against conventional mathematical modeling approaches, including nonlinear, linear, and convex formulations solved with established commercial software. \gls{mo} methods particularly excel in this domain due to their inherent capability to naturally handle the non-convex constraints that characterize power system operational equations, this is achieved by framing the optimization problem directly on a smooth manifold. Consequently, \gls{mo} offers a more robust framework for addressing these intricate challenges compared to heuristic evolutionary algorithms, convex relaxation techniques, or linearization strategies. The \gls{epds} case studies within this work specifically underscore \gls{mo}'s key strengths: its effective exploitation of the system's geometric properties and its capacity to manage the high-dimensional variable space typical of power flow problems, ultimately leading to accurate and efficient solutions.

The results from case studies using 14-bus, 33-bus, and 69-bus \gls{epds} models, and 3-bus and 4-bus \gls{epts} models, consistently highlight \gls{mo}'s effectiveness. The Manopt toolbox has demonstrated performance comparable to or even surpassing some commercial solvers in terms of execution time and accuracy for power flow problems in \gls{epds}. This indicates that \gls{mo} methods do not inherently lead to increased processing times and can maintain competitiveness in computational efficiency while delivering very high accuracy. Furthermore, \gls{mo} has shown robust performance across different system sizes, effectively managing the increasing complexity and dimensionality of modern power systems, thus positioning it as a robust solution for large-scale deployments.

Despite its promise, integrating \gls{mo} into existing power system infrastructures presents challenges related to implementation complexity, algorithmic design and scalability, in addition to future concerns such as planning for expansion. However, this paper has demonstrated that the practical implementation can be less intricate than initially perceived, especially with tools like Manopt and approaches like the \gls{bfs}. Moreover, \gls{mo}'s potential to incorporate other planning strategies, such as reconfiguration through Master-Slave configurations, further expands its applicability in power system optimization. For this, the Slave must have very low computational time and very high precision, and this is exactly what was proven in this article.

To fully harness \gls{mo}'s transformative potential in power systems optimization, continued research and development are essential. Addressing the remaining challenges through collaborative efforts from researchers and industry professionals will pave the way for more efficient, resilient, and scalable \gls{epds} and \gls{epts}.

\end{document}